\documentclass[11pt]{amsart}

\usepackage{amsmath}
\usepackage{amscd}
\usepackage{amssymb}
\usepackage{latexsym}
\usepackage{stmaryrd} 
\usepackage{mathbbol} 

\usepackage[arrow,curve,matrix,tips,2cell]{xy}
  \SelectTips{eu}{10} \UseTips
  \UseAllTwocells
\newtheorem{theorem}{Theorem}[section]
\newtheorem*{theorem*}{Theorem}
\newtheorem{lemma}[theorem]{Lemma}
\newtheorem*{lemma*}{Lemma}
\newtheorem{corollary}[theorem]{Corollary}
\newtheorem{proposition}[theorem]{Proposition}

\newtheorem{remark}[theorem]{Remark}
\newtheorem{definition}[theorem]{Definition}

\newcommand{\bgl}{\begin{equation}} 
\newcommand{\egl}{\end{equation}}
\newcommand{\bgloz}{\begin{equation*}} 
\newcommand{\egloz}{\end{equation*}}
\newcommand{\bgln}{\begin{eqnarray}} 
\newcommand{\egln}{\end{eqnarray}}
\newcommand{\bglnoz}{\begin{eqnarray*}} 
\newcommand{\eglnoz}{\end{eqnarray*}}
\newcommand{\btheo}{\begin{theorem}}
\newcommand{\etheo}{\end{theorem}}
\newcommand{\btheooz}{\begin{theorem*}}
\newcommand{\etheooz}{\end{theorem*}}
\newcommand{\blemma}{\begin{lemma}}
\newcommand{\elemma}{\end{lemma}}
\newcommand{\blemmaoz}{\begin{lemma*}}
\newcommand{\elemmaoz}{\end{lemma*}}
\newcommand{\bproof}{\begin{proof}}
\newcommand{\eproof}{\end{proof}}
\newcommand{\bbew}{\begin{beweis}}
\newcommand{\ebew}{\end{beweis}}
\newcommand{\bremark}{\begin{remark}\em}
\newcommand{\eremark}{\end{remark}}
\newcommand{\bdefin}{\begin{definition}}
\newcommand{\edefin}{\end{definition}}
\newcommand{\bprop}{\begin{proposition}}
\newcommand{\eprop}{\end{proposition}}
\newcommand{\bcor}{\begin{corollary}}
\newcommand{\ecor}{\end{corollary}}
\newcommand{\bfa}{\begin{cases}} 
\newcommand{\efa}{\end{cases}}
\newcommand{\cE}{\mathcal E}

\newcommand{\cI}{\mathcal I}
\newcommand{\cJ}{\mathcal J}

\newcommand{\cL}{\mathcal L}

\newcommand{\cV}{\mathcal V}
\newcommand{\cW}{\mathcal W}

\def\Cz{\mathbb{C}}

\def\Nz{\mathbb{N}}

\def\Zz{\mathbb{Z}}

\def\1z{\mathbb{1}}
\newcommand{\fT}{\mathfrak T}
\newcommand{\an}[1]{``#1''} 
\newcommand{\ti}{\tilde}

\newcommand{\lori}{\longrightarrow}

\newcommand{\larr}{\leftarrow}
\newcommand{\ma}{\mapsto} 
\newcommand{\mafr}{\mapsfrom} 
\newcommand\onto{\twoheadrightarrow} 
\newcommand\into{\hookrightarrow} 
\newcommand{\Rarr}{\Rightarrow} 
\newcommand{\LRarr}{\Leftrightarrow} 

\newcommand{\ve}{\varepsilon}

\def\SEMI{\mbox{$\times\kern-2pt\vrule height5pt width.6pt \kern3pt $}}

\newcommand{\halb}{\tfrac{1}{2}}

\newcommand{\End}{{\rm End}\,}
\newcommand{\Aut}{{\rm Aut}\,}

\newcommand{\id}{{\rm id}}

\newcommand{\stalg}{{}^*\text{-alg}}

\renewcommand{\ker}{{\rm ker}\,}

\newcommand{\reg}{^\times} 
\newcommand{\pos}{_{>0}} 
\newcommand{\lspan}{{\rm span}} 
\newcommand{\clspan}{\overline{\lspan}} 
\newcommand{\abs}[1]{\lvert#1\rvert} 
\newcommand{\norm}[1]{\left\|#1\right\|} 
\newcommand{\defeq}{\mathrel{:=}} 

\newcommand{\dop}{\text{: }} 
\newcommand{\falls}{\text{ if }} 
\newcommand{\sonst}{\text{ else}} 
\newcommand{\fa}{\text{ for all }} 
\newcommand{\e}[1]{e_{\left[#1\right]}} 
\newcommand{\E}[1]{E_{\left[#1\right]}} 
\newcommand{\rte}{{\rtimes}^{e}} 
\newcommand{\rta}{{\rtimes}^{a}} 
\newcommand{\un}{^{(\cup)}} 
\newcommand{\dsupp}{\text{d-supp}}
\newcommand{\Isom}{{\rm Isom}\,}
\newcommand{\Proj}{{\rm Proj}\,}

\newcommand{\lge}{\left\{} 
\newcommand{\rge}{\right\}} 
\newcommand{\lru}{\left(} 
\newcommand{\rru}{\right)} 
\newcommand{\lsp}{\left\langle} 
\newcommand{\rsp}{\right\rangle} 
\newcommand{\rukl}[1]{\lru #1 \rru} 
\newcommand{\gekl}[1]{\lge #1 \rge} 
\newcommand{\spkl}[1]{\lsp #1 \rsp} 
\newcommand{\menge}[2]{\gekl{ #1 \dop #2 }} 
\begin{document}

\title{Semigroup C*-algebras and amenability of semigroups}

\author{Xin Li}



\thanks{\scriptsize{Research supported by the Deutsche Forschungsgemeinschaft (SFB 878) and by the ERC through AdG 267079.}}

\begin{abstract}
We construct reduced and full semigroup C*-algebras for left cancellative semigroups. Our new construction covers particular cases already considered by A. Nica and also Toeplitz algebras attached to rings of integers in number fields due to J. Cuntz.

Moreover, we show how (left) amenability of semigroups can be expressed in terms of these semigroup C*-algebras in analogy to the group case.
\end{abstract}

\maketitle


\setlength{\parindent}{0pt} \setlength{\parskip}{0.5cm}

\section{Introduction}

The construction of group C*-algebras provides examples of C*-algebras which are both interesting and challenging to study. If we restrict our discussion to discrete groups, then we could say that the idea behind the construction is to implement the algebraic structure of a given group in a concrete or abstract C*-algebra in terms of unitaries. It then turns out that the group and its group C*-algebra(s) are closely related in various ways, for instance with respect to representation theory or in the context of amenability.

Given the success and the importance of the construction of group C*-algebras, a very natural question is whether we can start with algebraic structures that are even more basic than groups, namely semigroups. And indeed, this question has been addressed by various authors. The start was made by L. Coburn who studied the C*-algebra of the additive semigroup of the natural numbers (see \cite{Co1} and \cite{Co2}). Then, just to mention some examples, a number of authors like L. Coburn, R. G. Douglas, R. Howe, D. G. Schaeffer and I. M. Singer studied C*-algebras of particular Toeplitz operators in \cite{Co-Do}, \cite{C-D-S-S}, \cite{Dou} and \cite{Do-Ho}. The original motivation came from index theory and related K-theoretic questions. Later on, G. Murphy further generalized this construction, first to positive cones in ordered abelian groups in \cite{Mur1}, then to arbitrary left cancellative semigroups in \cite{Mur2} and \cite{Mur3}. The basic idea behind the constructions mentioned so far is to replace unitary representations in the group case by isometric representations for left cancellative semigroups. However, it turns out that the full semigroup C*-algebras introduced by G. Murphy are very complicated and not suited for studying amenability. For instance, the full semigroup C*-algebra of $\Nz \times \Nz$ in the sense of G. Murphy is not nuclear (see \cite{Mur4}, Theorem~6.2).

Apart from these constructions, A. Nica has introduced a different construction of semigroup C*-algebras for positive cones in quasi-lattice ordered groups (see \cite{Ni} and also \cite{La-Rae}). His construction has the advantage that it leads to much more tractable C*-algebras than the construction introduced by G. Murphy, so that A. Nica was able to study amenability questions using his new construction. The main difference between A. Nica's construction and the former ones is that A. Nica takes the right ideal structure of the semigroups into account in his construction, although in a rather implicit way.

Another source of inspiration is provided by so-called ring C*-algebras (see \cite{Cun}, \cite{Cu-Li1}, \cite{Cu-Li2} and \cite{Li}). Namely, the author realized during his recent work \cite{Li} that there are strong parallels between the construction of ring C*-algebras and semigroup C*-algebras. The restriction A. Nica puts on his semigroups by only considering positive cones in quasi-lattice ordered groups would correspond in the ring case to considering rings for which every ideal is principal. This observation indicates that the ideal structure (of the ring or semigroup) should play an important role in more general constructions. This idea has been worked out in the case of rings in \cite{Li}. Moreover, it was explained in Appendix A.2 of \cite{Li} how the analogous idea leads to a generalization of A. Nica's construction to arbitrary left cancellative semigroups.

Independently from this construction of semigroup C*-algebras, J. Cuntz has modified the construction of ring C*-algebras from \cite{Cu-Li1} and \cite{Cu-Li2} and has introduced so-called Toeplitz algebras for certain rings from algebraic number theory (rings of integers in number fields). The motivation was to improve the functorial properties of ring C*-algebras. And again, the crucial idea behind the construction is to make use of the ideal structure of the rings of interest. This first step was due to J. Cuntz (before the work \cite{C-D-L}), and he presented these ideas and the results on functoriality in a talk at the \an{Workshop on C*-algebras} in Nottingham which took place in September 2010.

As a next step, J. Cuntz, C. Deninger and M. Laca study these Toeplitz algebras in \cite{C-D-L} and they show that the Toeplitz algebra of the ring of integers in a number field can be identified via a canonical representation with the reduced semigroup C*-algebra of the $ax+b$-semigroup over the ring. This indicates that there is a strong connection between these Toeplitz algebras and semigroup C*-algebras.

And indeed, it turns out that if we apply the construction of full semigroup C*-algebras in \cite{Li} to the $ax+b$-semigroups over rings of integers, then we arrive at universal C*-algebras which are canonically isomorphic to these Toeplitz algebras. As pointed out in \cite{C-D-L}, the most interesting examples arise from rings which do not have the property that every ideal is principal (i.e. the class number of the number field is strictly bigger than $1$). For these rings or rather the corresponding $ax+b$-semigroups, it is not possible to apply A. Nica's construction. This explains the need for a generalization of A. Nica's work.

So, to summarize, the motivation behind our construction of semigroup C*-algebras is twofold: On the one hand, we would like to provide a general framework for A. Nica's constructions as well as the Toeplitz algebras due to J. Cuntz so that these constructions can be naturally thought of as particular cases of our general construction (this is explained in \S~\ref{constructions}). On the other hand, we would like to obtain constructions which are more tractable than those of G. Murphy and which allow us to characterize amenability of semigroups very much in the same spirit as in the group case (see \S~\ref{amenability}). To establish this connection with amenability, we first have to modify our construction of full semigroup C*-algebras in the case of subsemigroups of groups (see \S~\ref{variant}).

Of course, there are not only C*-algebras associated with groups, but also C*-algebras attached to dynamical systems. So another question would be whether we can also construct C*-algebras for semigroup actions by automorphisms. We only touch upon this question in \S~\ref{cropro-auto}.


I would like to thank J. Cuntz for interesting and helpful discussions and for providing access to the preprint \cite{C-D-L} due to him, C. Deninger and M. Laca. I also thank M. Norling who has pointed me towards a missing relation in the definition of full semigroup C*-algebras for subsemigroups of groups. This has led me to the modified construction introduced in \S~\ref{variant}.

\section{Constructions}
\label{constructions}

\subsection{Semigroup C*-algebras}
By a semigroup, we mean a set $P$ equipped with a binary operation $P \times P \to P; (p,q) \ma pq$ which is associative, i.e. $(p_1 p_2)p_3 = p_1(p_2 p_3)$. We always assume that our semigroup has a unit element, i.e. there exists $e \in P$ such that $ep=pe=p$ for all $p \in P$. All semigroup homomorphisms shall preserve unit elements. We only consider discrete semigroups. A semigroup $P$ is called left cancellative if for every $p$, $x$ and $y$ in $P$, $px=py$ implies $x=y$.

As mentioned in the introduction, the basic idea behind the construction of semigroup C*-algebras is to represent semigroup elements by isometries. This means that if we let $\Isom$ be the semigroup of the necessarily unital semigroup C*-algebra associated with the semigroup $P$, then we would like to have a semigroup homomorphism $P \to \Isom$. This requirement explains why we restrict our discussion to left cancellative semigroups: Since $\Isom$ is always a left cancellative semigroup, this homomorphism $P \to \Isom$ can only be faithful if $P$ itself is left cancellative.

Given a left cancellative semigroup $P$, we can construct its left regular representation as follows:

Let $\ell^2(P)$ be the Hilbert space of square summable complex-valued functions on $P$. $\ell^2(P)$ comes with the canonical orthonormal basis $\menge{\varepsilon_x}{x \in P}$ given by $\varepsilon_x(y) = \delta_{x,y}$ where $\delta_{x,y}$ is $1$ if $x=y$ and $0$ if $x \neq y$. Let us define for every $p \in P$ an isometry $V_p$ by setting $V_p \varepsilon_x = \varepsilon_{px}$. Here we have made use of our assumption that our semigroup $P$ is left cancellative. It ensures that the assignment $\varepsilon_x \ma \varepsilon_{px}$ indeed extends to an isometry. Now the reduced semigroup C*-algebra of $P$ is simply given as the sub-C*-algebra of $\cL(\ell^2(P))$ generated by these isometries $\menge{V_p}{p \in P}$. We denote this concrete C*-algebra by $C^*_r(P)$, i.e. we set
\bdefin
$C^*_r(P) \defeq C^* \rukl{\menge{V_p}{p \in P}} \subseteq \cL(\ell^2(P))$.
\edefin
So $C^*_r(P)$ is really a very natural object: It is the C*-algebra generated by the left regular representation of $P$. This C*-algebra $C^*_r(P)$ is called the reduced semigroup C*-algebra of $P$ in analogy to the group case. But we remark that this C*-algebra is also called the Toeplitz algebra of $P$ by various authors.

We now turn to the construction of full semigroup C*-algebras. As explained in the introduction, we will make use of right ideals of our semigroups to construct full semigroup C*-algebras. So we first have to choose a family of right ideals.

Given a semigroup $P$, every semigroup element $p \in P$ gives rise to the map $P \to P; x \ma px$. It is simply given by left multiplication with $p$. Given a subset $X$ of $P$ and an element $p \in P$, we set 
\bgl
\label{notations}
  pX \defeq \menge{px}{x \in X} \text{ and } p^{-1}X \defeq \menge{y \in P}{py \in X}.
\egl
In other words, $pX$ is the image and $p^{-1}X$ is the pre-image of $X$ under left multiplication with $p$. A subset $X$ of $P$ is called a right ideal if it is closed under right multiplication with arbitrary semigroup elements, i.e. if for every $x \in X$ and $p \in P$, the product $xp$ always lies in $X$. 

The semigroup $P$ is left cancellative if and only if for every $p \in P$, left multiplication with $p$ defines an injective map. For the rest of this section, let $P$ always be a left cancellative semigroup.

Let $\cJ$ be the smallest family of right ideals of $P$ containing $P$ and $\emptyset$, i.e. 
\bgl
\label{J-cond-P-0}
  P \in \cJ, \emptyset \in \cJ,
\egl
and closed under left multiplication, taking pre-images under left multiplication,
\bgl
\label{J-cond-lm}
  X \in \cJ, p \in P \Rarr pX, p^{-1}X \in \cJ,
\egl
as well as finite intersections,
\bgl
\label{J-cond-int}
  X, Y \in \cJ \Rarr X \cap Y \in \cJ.
\egl
It is not difficult to find out how right ideals in $\cJ$ typically look like. Actually, it follows directly from the definitions that
\bgl
\label{J-explicit}
  \cJ = \menge{\bigcap_{j=1}^N (q_{j,1})^{-1} p_{j,1} \dotsm (q_{j,n_j})^{-1} p_{j,n_j} P}{N, n_j \in \Zz \pos; p_{j,k}, q_{j,k} \in P} \cup \gekl{\emptyset}.
\egl
The elements in $\cJ$ are called constructible right ideals. If we want to keep track of the semigroup, we write $\cJ_P$ for the family of constructible right ideals of the semigroup $P$. We will see in \eqref{J-without-int} that it is not necessary to ask for \eqref{J-cond-int}.

With the help of this family of right ideals, we can now construct the full semigroup C*-algebra of $P$. The idea is to ask for a projection-valued spectral measure, defined for elements in the family $\cJ$ and taking values in projections in our C*-algebra.

\bdefin
\label{full-semigp-C}
The full semigroup C*-algebra of $P$ is the universal C*-algebra generated by isometries $\menge{v_p}{p \in P}$ and projections $\menge{e_X}{X \in \cJ}$ satisfying the following relations:
\bglnoz
  && \text{I.(i) } v_{pq} = v_p v_q \ \ \ \text{I.(ii) } v_p e_X v_p^* = e_{pX} \\
  && \text{II.(i) } e_P = 1 \ \ \ \text{II.(ii) } e_{\emptyset} = 0 \ \ \ \text{II.(iii) } e_{X \cap Y} = e_X \cdot e_Y
\eglnoz
for all $p$, $q$ in $P$ and $X$, $Y$ in $\cJ$.

We denote this universal C*-algebra by $C^*(P)$, i.e.
\bgloz
  C^*(P) \defeq C^* \rukl{\menge{v_p}{p \in P} \cup \menge{e_X}{X \in \cJ} \vline 
  \begin{array}{c}
  v_p \text{ are isometries }
  \\
  \text{and } e_X \text{ are projections }
  \\
  \text{satisfying I and II.}
  \end{array}
  }
\egloz
\edefin
One remark about notation: For the sake of readability, we sometimes write $\e{X}$ for $e_X$ in case the expression in the index gets very long.

Of course, the question is: Where do all these relations come from? The idea is that we can think of $C^*(P)$ as a universal model of the reduced semigroup C*-algebra $C^*_r(P)$. To make this precise, let us again consider concrete operators on $\ell^2(P)$. We have already defined the isometries $V_p$ for $p \in P$. For every subset $X$ of $P$, let $E_X$ be the orthogonal projection onto $\ell^2(X) \subseteq \ell^2(P)$. In other words, let $\1z_X$ be the characteristic function of $X$ defined on $P$, i.e. $\1z_X(p) = 1$ if $p \in X$ and $\1z_X(p) = 0$ if $p \notin X$. Then $\1z_X$ is an element of $\ell^\infty(P)$ which is mapped to $E_X$ under the canonical representation of $\ell^\infty(P)$ as multiplication operators on $\ell^2(P)$. As with the projections $e_X$, we will sometimes write $\E{X}$ for $E_X$ if the subscript becomes very long. It is now easy to check that the two families $\menge{V_p}{p \in P}$ and $\menge{E_X}{X \in \cJ}$ satisfy relations I and II (with $V_p$ in place of $v_p$ and $E_X$ in place of $e_X$). This explains the origin of these relations. At the same time, we obtain by universal property of $C^*(P)$ a non-zero homomorphism $\lambda: C^*(P) \to \cL(\ell^2(P))$ sending $v_p$ to $V_p$ and $e_X$ to $E_X$ for every $p \in P$ and $X \in \cJ$. This homomorphism is called the left regular representation of $C^*(P)$. In particular, we see that $C^*(P)$ is not the zero C*-algebra. We will see later on (compare \eqref{cd}) that the image of $\lambda$ is actually the reduced semigroup C*-algebra $C^*_r(P)$.

\bremark
\label{cJ'}
Actually, the requirement that $\cJ$ should be closed under taking pre-images under left multiplications is not needed in the construction, and it does not appear in the first version of semigroup C*-algebras in \cite{Li}, Appendix~A.2. The original reason why we added this extra requirement is that we wanted our construction of full semigroup C*-algebras to include the construction of Toeplitz algebras for rings of integers in number fields by J. Cuntz. However, for such semigroups, it is not necessary to consider pre-images in the following sense: Let $\cJ'$ denote the family of right ideals defined in the same way as $\cJ$ but without the property that $\cJ'$ is closed under pre-images under left multiplication. For the $ax+b$-semigroups over rings of integers, it turns out that it does not matter whether we take $\cJ$ or $\cJ'$ in Definition~\ref{full-semigp-C} because the resulting C*-algebras are canonically isomorphic. But for general semigroups, it is more convenient to work with $\cJ$ as we will see.
\eremark

Let us also discuss a useful modification of these full semigroup C*-algebras. We first reformulate relation II.(iii): We have canonical lattice structures on the set of right ideals of $P$ (let $X \wedge Y = X \cap Y$ and $X \vee Y = X \cup Y$ for right ideals $X$ and $Y$) and on the set of commuting projections in a C*-algebra (let $e \wedge f = ef$ and $e \vee f = e+f-e \wedge f$ for commuting projections $e$ and $f$). So relation II.(iii) simply tells us that the projections $\menge{e_X}{X \in \cJ}$ commute and that the assignment $\cJ \ni X \ma e_X \in \Proj(C^*(P))$ is $\wedge$-compatible. Given this interpretation, an obvious question is whether we can modify our construction so that the analogous assignment becomes $\vee$-compatible as well. This is indeed possible. The first step is to enlarge the family $\cJ$ so that it is closed under finite unions as well. Let $\cJ\un$ be the smallest family of right ideals of $P$ satisfying the conditions \eqref{J-cond-P-0} -- \eqref{J-cond-int} and the extra condition
\bgl
\label{J-cond-union}
  X,Y \in \cJ\un \Rarr X \cup Y \in \cJ\un.
\egl
Again, it follows from our definition that
\bgl
\label{J-union-explicit}
  \cJ\un = \menge{\bigcup_{i=1}^M \bigcap_{j=1}^N (q_{j,1}^{(i)})^{-1} p_{j,1}^{(i)} \dotsm (q_{j,n_j}^{(i)})^{-1} p_{j,n_j}^{(i)} P}
  {M,N, n_j \in \Zz \pos; p_{j,k}^{(i)}, q_{j,k}^{(i)} \in P} \cup \gekl{\emptyset}.
\egl

We can now modify Definition~\ref{full-semigp-C} by replacing $\cJ$ by $\cJ\un$ and adding to the relations the extra relation $e_{X \cup Y} = e_X + e_Y - e_{X \cap Y}$ for all $X,Y \in \cJ\un$. The corresponding universal C*-algebra is then denoted by ${C^*}\un(P)$.
\bdefin
\label{full-semigp-C-union}
\bgloz
  {C^*}\un(P) \defeq C^* \rukl{\menge{v_p}{p \in P} \cup \menge{e_X}{X \in \cJ\un} \vline 
  \begin{array}{c}
  v_p \text{ are isometries }
  \\
  \text{and } e_X \text{ are projections }
  \\
  \text{satisfying I and } \text{II}\un
  \end{array}
  }
\egloz
with the relations
\bglnoz
  && \text{I.(i) } v_{pq} = v_p v_q \ \ \ \text{I.(ii) } v_p e_X v_p^* = e_{pX} \\
  && \text{II}\un.\text{(i) } e_P = 1 \ \ \ \text{II}\un.\text{(ii) } e_{\emptyset} = 0 \\
  && \text{II}\un.\text{(iii) } e_{X \cap Y} = e_X \cdot e_Y \ \ \ \text{II}\un.\text{(iv) } e_{X \cup Y} = e_X + e_Y - e_{X \cap Y}.
\eglnoz
\edefin
It is immediate from our definitions that ${C^*}\un(P)$ is a quotient of $C^*(P)$, or in other words, that we always have a canonical homomorphism $\pi\un: C^*(P) \to {C^*}\un(P)$ sending $C^*(P) \ni v_p$ to $v_p \in {C^*}\un(P)$ and $C^*(P) \ni e_X$ to $e_X \in {C^*}\un(P)$ for all $p \in P$ and $X \in \cJ \subseteq \cJ\un$. Relation II$\un$.(iv) implies that $\pi\un$ is always surjective.

As for the relations defining $C^*(P)$, it is immediate that the relations I and II$\un$ (with $V_p$ in place of $v_p$ and $E_X$ in place of $e_X$) are satisfied by the concrete operators $\menge{V_p}{p \in P}$ and $\menge{E_X}{X \in \cJ\un}$ on $\ell^2(P)$ ($E_X$ is the orthogonal projection onto $\ell^2(X) \subseteq \ell^2(P)$ as above). So we again obtain by universal property of ${C^*}\un(P)$ a non-zero homomorphism $\lambda\un: {C^*}\un(P) \to \cL(\ell^2(P))$ sending $v_p$ to $V_p$ and $e_X$ to $E_X$ for every $p \in P$ and $X \in \cJ\un$. This again implies that ${C^*}\un(P)$ is not the zero C*-algebra. Moreover, we obtain by construction a commutative diagram
\bgl
\xymatrix{
C^*(P) \ar[d]_{\pi\un} \ar[dr]^{\lambda}
&
\\
{C^*}\un(P) \ar[r]_{\lambda\un} & \cL(\ell^2(P))
}
\label{cd-to-cL}
\egl

\subsection{Semigroup crossed products by automorphisms}
\label{cropro-auto}

At this point, we also introduce semigroup crossed products by automorphisms. Let $P$ be a left cancellative semigroup and $D$ a unital C*-algebra. Moreover, let $\alpha: P \to \Aut(A)$ be a semigroup homomorphism.

We then define the full semigroup crossed product of $A$ by $P$ with respect to $\alpha$ as the (up to isomorphism unique) unital C*-algebra $A \rta_\alpha P$ which comes with two unital homomorphisms $\iota_A: A \to A \rta_\alpha P$ and $\iota_P: C^*(P) \to A \rta_\alpha P$ satisfying
\bgloz
  \iota_A(\alpha_p(a)) \iota_P(v_p) = \iota_P(v_p) \iota_A(a) \fa a \in A, p \in P
\egloz
such that the following universal property is fulfilled:

Whenever $T$ is a unital C*-algebra and $\varphi_A: A \to T$, $\varphi_P: C^*(P) \to T$ are unital homomorphisms satisfying the covariance relation
\bgl
\label{cov-rel}
  \varphi_A(\alpha_p(a)) \varphi_P(v_p) = \varphi_P(v_p) \varphi_A(a) \fa a \in A, p \in P,
\egl
there is a unique homomorphism $\varphi_A \rtimes \varphi_P: A \rta_\alpha P \to T$ with
\bgloz
  (\varphi_A \rtimes \varphi_P) \circ \iota_A = \varphi_A \text{ and } (\varphi_A \rtimes \varphi_P) \circ \iota_P = \varphi_P.
\egloz
We could also use ${C^*}\un(P)$ instead of $C^*(P)$ in the construction of the semigroup crossed product by automorphisms, and the result would be another C*-algebra, say $A \rtimes^{a,(\cup)}_\alpha P$, with the corresponding universal property. We will see in Lemma~\ref{rta-exists} that these universal C*-algebras really exist. By construction, we have a canonical homomorphism $\pi\un_{(A,P,\alpha)}: A \rta_\alpha P \to A \rtimes^{a,(\cup)}_\alpha P$. This homomorphism is surjective as the canonical homomorphism $\pi\un: C^*(P) \to {C^*}\un(P)$ is surjective. Of course, if $\text{tr}: P \to \Aut(\Cz)$ denotes the trivial action, then
\bgloz
  C^*(P) \cong \Cz \rtimes_{\text{tr}} P \text{, } {C^*}\un(P) \cong \Cz \rtimes\un_{\text{tr}} P,
\egloz
and under these canonical identifications, $\pi\un_{(\Cz,P,\text{tr})}$ becomes the canonical homomorphism $\pi\un: C^*(P) \to {C^*}\un(P)$.

We remark that there is a different notion of semigroup crossed products by endomorphisms which is for instance explained in \cite{La}, \cite{La-Rae}, \S~2 or in \cite{Li}, Appendix~A.1. We denote semigroup crossed products by endomorphisms by $\rte$ to distinguish them from our construction. We will see that there is a close relationship between these two sorts of semigroup crossed products.

G. Murphy has already introduced semigroup crossed products by automorphisms in \cite{Mur2} and \cite{Mur3}. However, as in the case of semigroup C*-algebras, G. Murphy's construction leads to very complicated C*-algebras which are not tractable even in very simple cases. But G. Murphy has also constructed concrete representations, and these can be used to define reduced semigroup crossed products by automorphisms: Take a faithful representation of $D$ on a Hilbert space $H$, say $i: A \to \cL(H)$. Form the tensor product $H \otimes \ell^2(P)$. Then define for every $a$ in $A$ a bounded operator by the formula $\eta \otimes \varepsilon_x \ma i(\alpha_x^{-1}(a)) (\eta) \otimes \varepsilon_x$ for every $\eta \in H$ and $x \in P$. It is straightforward to check that these operators give rise to a homomorphism $i_A: A \to \cL(H \otimes \ell^2(P))$ and that $i_A$ and $i_P \defeq \id_H \otimes \lambda: C^*(P) \to \cL(H \otimes \ell^2(P))$ satisfy the covariance relation \eqref{cov-rel}. Thus we obtain by universal property of $A \rta_\alpha P$ a homomorphism $\lambda_{(A,P,\alpha)} \defeq i_A \rtimes i_P: A \rta_\alpha P \to \cL(H \otimes \ell^2(P))$. We set $A \rta_{\alpha,r} P \defeq \lambda_{(A,P,\alpha)}(A \rta_\alpha P)$ and call this algebra the reduced semigroup crossed product of $A$ by $P$ with respect to $\alpha$. Using the same faithful representation $i$ of $A$, the induced homomorphism $i_A: A \to \cL(H \otimes \ell^2(P))$ and the homomorphism $\id_H \otimes \lambda\un: {C^*}\un(P) \to \cL(H \otimes \ell^2(P))$, we can also construct a homomorphism $\lambda_{(A,P,\alpha)}\un: A \rtimes^{a,(\cup)}_\alpha P \to \cL(H \otimes \ell^2(P))$. Again, by universal property of $A \rta_\alpha P$, $\lambda_{(A,P,\alpha)} = \lambda_{(A,P,\alpha)}\un \circ \pi\un_{(A,P,\alpha)}$, so there is no difference between $A \rtimes^{a,(\cup)}_{\alpha,r} P \defeq \lambda_{(A,P,\alpha)}\un (A \rtimes^{a,(\cup)}_\alpha P)$ and $A \rta_{\alpha,r} P$.

\bremark
Of course, we can consider right cancellative semigroups instead of left cancellative ones. Replacing left multiplication by right multiplication and right ideals by left ideals, we obtain analogous constructions. Alternatively, given a right cancellative semigroup $P$, we can go over to the opposite semigroup $P^{\text{op}}$ consisting of the same underlying set $P$ equipped with a new binary operation $\bullet$ given by $p \bullet q \defeq qp$. It is immediate that $P^{\text{op}}$ is left cancellative and our constructions apply.

With the obvious modifications, our analysis of C*-algebras associated with left cancellative semigroups (which is going to come) carries over to right cancellative semigroups.
\eremark

\subsection{Direct consequences of the definitions}
First of all, each of the C*-algebras $C^*(P)$ and ${C^*}\un(P)$ contains a distinguished sub-C*-algebra, namely the one generated by the projections $\menge{e_X}{X \in \cJ} \text{ or } \menge{e_X}{X \in \cJ\un}$. Let us denote these sub-C*-algebras by $D(P)$ and $D\un(P)$, i.e. 
\bglnoz
  && D(P) \defeq C^*(\menge{e_X}{X \in \cJ}) \subseteq C^*(P) \\
  && D\un(P) \defeq C^*(\menge{e_X}{X \in \cJ\un}) \subseteq {C^*}\un(P).
\eglnoz
We first observe that 
\bgl
\label{pi(D)}
\pi\un(D(P)) = D\un(P).
\egl
The inclusion \an{$\subseteq$} is clear as $\cJ \subseteq \cJ\un$, and the reverse inclusion \an{$\supseteq$} follows immediately from relation II$\un$.(iv) and the concrete description of $\cJ\un$ in \eqref{J-union-explicit}.

Moreover, we have the following
\blemma
The families $\menge{e_X}{X \in \cJ}$ and $\menge{e_X}{X \in \cJ\un}$ consist of commuting projections and are multiplicatively closed.
\elemma
\bproof
This follows immediately from relation II.(iii) and II$\un$.(iii), respectively.
\eproof
\bcor
$D(P)$ and $D\un(P)$ are commutative C*-algebras.

Moreover, $D(P)=\clspan(\menge{e_X}{X \in \cJ})$ and $D\un(P)=\clspan(\menge{e_X}{X \in \cJ\un})$.
\ecor
Furthermore, as another consequence of the definitions, we derive
\blemma
\label{v*_v}
For every $p \in P$ and $X \in \cJ$ ($X \in \cJ\un$), we have $v_p^* e_X v_p = e_{p^{-1} X}$ in $C^*(P)$ (${C^*}\un(P)$).
\elemma
\bproof
The proof is the same for $C^*(P)$ and ${C^*}\un(P)$. Take $p \in P$ and $X \in \cJ$ ($X \in \cJ\un$). We then have $v_p^* e_X v_p = v_p^* e_X v_p v_p^* v_p = v_p^* e_X e_{pP} v_p = v_p^* e_{X \cap pP} v_p = v_p^* e_{p(p^{-1} X)} v_p = v_p^* v_p e_{p^{-1} X} v_p^* v_p = e_{p^{-1} X}$.
\eproof
\bcor
For every $p \in P$, conjugation by $v_p^* \in C^*(P)$ ($v_p^* \in {C^*}\un(P)$) induces a homomorphism on $D(P)$ ($D\un(P)$).
\ecor
\bproof
This is a direct consequence of the previous lemma.
\eproof
From Lemma~\ref{v*_v} and the description of $\cJ$ given in \eqref{J-explicit}, we immediately deduce
\bcor
\label{v_p-generate}
$C^*(P)$ is generated as a C*-algebra by the isometries $\menge{v_p}{p \in P}$.
\ecor
We also obtain the analoguous statement for ${C^*}\un(P)$:
\bcor
\label{v_p-generate-union}
${C^*}\un(P)$ is generated as a C*-algebra by $\menge{v_p}{p \in P}$.
\ecor
\bproof
This either follows analogously from Lemma~\ref{v*_v} for ${C^*}\un(P)$ and the explicit description of $\cJ\un$ in \eqref{J-union-explicit} or with the help of the last corollary and the surjection $\pi\un: C^*(P) \to {C^*}\un(P)$.
\eproof
Now, it follows from Corollary~\ref{v_p-generate} that the image of the left regular representation $\lambda: C^*(P) \to \cL(\ell^2(P))$ is precisely the reduced semigroup C*-algebra $C^*_r(P)$. This means that we can rewrite the commutative triangle \eqref{cd-to-cL} more accurately as follows:
\bgl
\xymatrix{
C^*(P) \ar[d]_{\pi\un} \ar[dr]^{\lambda}
&
\\
{C^*}\un(P) \ar[r]_{\lambda\un} & C^*_r(P)
}
\label{cd}
\egl
As we did before for the full semigroup C*-algebras, we consider a canonical sub-C*-algebra of $C^*_r(P)$:
\bdefin
$D_r(P) \defeq C^*(\menge{E_X}{X \in \cJ}) \subseteq \cL(\ell^2(P))$.
\edefin
Recall that $E_X$ is the orthogonal projection onto the subspace $\ell^2(X) \subseteq \ell^2(P)$.

It is immediately clear that $\lambda(D(P)) = D_r(P)$, so that $D_r(P)$ is a sub-C*-algebra of $C^*_r(P)$. $D_r(P)$ is obviously commutative and we have $D_r(P) = \clspan(\menge{E_X}{X \in \cJ})$ since $\menge{E_X}{X \in \cJ}$ is multiplicatively closed. Because of $\lambda(D(P)) = D_r(P)$, the commutative triangle \eqref{cd}, restricted to the distinguished commutative sub-C*-algebras, yields the commutative triangle
\bgl
\xymatrix{
D(P) \ar[d]_{\pi\un} \ar[dr]^{\lambda}
&
\\
D\un(P) \ar[r]_{\lambda\un} & D_r(P)
}
\label{cd-d}
\egl
Another direct consequence of our constructions is that we can alternatively describe our constructions as semigroup crossed products by endomorphisms. For the reader's convenience, we recall the notion of semigroup crossed products by endomorphisms. Let $P$ be a discrete semigroup and $D$ a unital C*-algebra. Further assume that $\tau: P \to \End(D)$ is a semigroup homomorphism from $P$ to the semigroup $\End(D)$ of (not necessarily unital) endomorphisms of $D$.
\bdefin
The semigroup crossed product $D \rte_{\tau} P$ is the up to canonical isomorphism unique unital C*-algebra which comes with a unital homomorphism $i_D: D \to D \rte_{\tau} P$ and a semigroup homomorphism $i_P: P \to \Isom(D \rte_{\tau} P)$ subject to the condition $i_P(p) i_D(d) i_P(p)^* = i_D(\tau_p(a))$ for all $p \in P$, $d \in D$ and satisfying the following universal property:

Whenever $T$ is a unital C*-algebra, $j_D: D \to T$ is a unital homomorphism and $j_P: P \to \Isom(T)$ is a semigroup homomorphism such that the covariance relation
\bgl
\label{cov-rel-endo}
  j_P(p) j_D(d) j_P(p)^* = j_D(\tau_p(d)) \fa p \in P, d \in D
\egl
is fulfilled, there is a unique homomorphism $j_D \rtimes j_P: D \rte_{\tau} P \to T$ with $(j_D \rtimes j_P) \circ i_D = j_D$ and $(j_D \rtimes j_P) \circ i_P = j_P$. Here $\Isom(D \rte_{\tau} P)$ and $\Isom(T)$ are the semigroups of isometries in $D \rte_{\tau} P$ and $T$, respectively.
\edefin
Existence of $D \rte_{\tau} P$ is shown in \cite{La-Rae}, \S~2; their condition (iii) is equivalent to uniqueness of $j_D \rtimes j_P$.

Now, in our situation, there are canonical actions (i.e. semigroup homomorphisms) $\tau: P \to \End(D(P))$ and $\tau\un: P \to \End(D\un(P))$ given by $P \ni p \ma v_p \sqcup v_p^*$. Conjugation by $v_p$ gives rise to a homomorphism of $C^*(P)$ because $v_p$ is an isometry, and $D(P)$ ($D\un(P)$) is invariant under these homomorphisms by relation I.(ii). When we form the corresponding semigroup crossed products by endomorphisms, we obtain
\blemma
\label{cropro-descriptions}
$C^*(P)$ is canonically isomorphic to $D(P) \rte_{\tau} P$, and ${C^*}\un(P)$ is canonically isomorphic to $D\un(P) \rte_{\tau\un} P$.
\elemma
\bproof
Using the universal property of $C^*(P)$ and $D(P) \rte_{\tau} P$, we can construct mutually inverse homomorphisms $C^*(P) \leftrightharpoons D(P) \rte_{\tau} P$. It is clear that the isometries $\menge{i_P(p)}{p \in P} \subseteq D(P) \rte_{\tau} P$ and the projections $\menge{i_{D(P)}(e_X)}{X \in \cJ} \subseteq D(P) \rte_{\tau} P$ satisfy relations I and II (in place of the $v_p$s and $e_X$s), so that there exists a homomorphism $C^*(P) \to D(P) \rte_{\tau} P$ sending $v_p$ to $i_P(p)$ and $e_X$ to $i_{D(P)}(e_X)$ for all $p \in P$ and $X \in \cJ$. Conversely, $C^*(P)$ together with the inclusion $D(P) \into C^*(P)$ and the semigroup homomorphism $P \ni p \ma v_p \in \Isom(C^*(P))$ satisfies the covariance relation~\eqref{cov-rel-endo} because of relation I.(ii). Hence there exists a homomorphism $D(P) \rte_{\tau} P \to C^*(P)$ sending $i_P(p)$ to $v_p$ and $i_{D(P)}(e_X)$ to $e_X$ for all $p \in P$ and $X \in \cJ$. By construction, these two homomorphisms are inverse to one another.

Similarly, a comparison of the universal properties yields a canonical identification ${C^*}\un(P) \cong D\un(P) \rte_{\tau\un} P$.
\eproof
More generally, we can also describe $D \rta_{\alpha} P$ and $D \rtimes^{a,(\cup)}_\alpha P$ as crossed products.
\blemma
\label{rta-exists}
$A \rta_{\alpha} P$ and $A \rtimes^{a,(\cup)}_{\alpha} P$ exist and are canonically isomorphic to $(A \otimes D(P)) \rte_{\alpha \otimes \tau} P$ and $(A \otimes D\un(P)) \rte_{\alpha \otimes \tau\un} P$, respectively.
\elemma
\bproof
By construction, $A \rta_{\alpha} P$ and $(A \otimes D(P)) \rte_{\alpha \otimes \tau} P$ have the same universal property. (Note that relation~\eqref{cov-rel} implies that $\iota_A(A)$ and $\iota_P(D(P))$ in $A \rta_{\alpha} P$ commute.) As $(A \otimes D(P)) \rte_{\alpha \otimes \tau} P$ exists by \cite{La-Rae}, Proposition~2.1, we have proven our assertions about $A \rta_{\alpha} P$. An analogous argument applies to $A \rtimes^{a,(\cup)}_{\alpha} P$.
\eproof
Another observation is that our constructions behave nicely with respect to direct products of semigroups.
\blemma
\label{C(product)}
Given two left cancellative semigroups $P$ and $Q$, there are canonical isomorphisms
\bglnoz
  && C^*(P \times Q) \cong C^*(P) \otimes_{\max} C^*(Q) \text{ given by } v_{(p,q)} \ma v_p \otimes v_q \\
  &\text{and}& C^*_r(P \times Q) \cong C^*_r(P) \otimes_{\min} C^*_r(Q) \text{ given by } V_{(p,q)} \ma V_p \otimes V_q.
\eglnoz
\elemma
\bproof
For the first identification, we just have to compare the universal properties of these C*-algebras. The second identification is given by conjugation by the unitary $\ell^2(P) \otimes \ell^2(Q) \to \ell^2(P \times Q); \varepsilon_x \otimes \varepsilon_y \ma \varepsilon_{(x,y)}$.
\eproof
\bremark
\label{C-un(product)1}
We can also identify ${C^*}\un(P \times Q)$ with ${C^*}\un(P) \otimes_{\max} {C^*}\un(Q)$ via $v_{(p,q)} \ma v_p \otimes v_q$. The problem is to show that there is a homomorphism $D\un(P \times Q) \to {C^*}\un(P) \otimes_{\max} {C^*}\un(Q)$ which sends for all $X \in \cJ_P$ and $Y \in \cJ_Q$ the projection $e_{X \times Y}$ to $e_X \otimes e_Y$. This has to be the case as we want that $v_{(p,q)}$ is sent to $v_p \otimes v_q$ for every $p \in P$ and $q \in Q$. Once we know that such a homomorphism $D\un(P \times Q) \to {C^*}\un(P) \otimes_{\max} {C^*}\un(Q)$ exists, we can easily construct, using Lemma~\ref{cropro-descriptions}, the desired homomorphism ${C^*}\un(P \times Q) \to {C^*}\un(P) \otimes_{\max} {C^*}\un(Q)$ satisfying $v_{(p,q)} \ma v_p \otimes v_q$. It is also easy to construct the inverse homomorphism ${C^*}\un(P \times Q) \larr {C^*}\un(P) \otimes_{\max} {C^*}\un(Q)$. It turns out that such a desired homomorphism $D\un(P \times Q) \to {C^*}\un(P) \otimes_{\max} {C^*}\un(Q)$ indeed exists (see Corollary~\ref{C-un(product)2}). But the proof will have to wait until we have studied in more detail the relationship between $D\un(P)$ and $D_r(P)$.
\eremark

\subsection{Examples}
\label{exp}
Of course, if $P$ happens to be a group, then our constructions coincide with the usual constructions of group C*-algebras or ordinary crossed products. To be more precise, if $P$ is a group, then the canonical homomorphism $\pi\un: C^*(P) \to {C^*}\un(P)$ is an isomorphism. Moreover, $C^*(P)$ and $C^*_r(P)$ can be canonically identified with the full and the reduced group C*-algebra of the group $P$. Analogously, for every unital C*-algebra $A$ and every (semi)group homomorphism $P \to \Aut(A)$, the canonical homomorphism $\pi\un_{(A,P,\alpha)}: A \rta_\alpha P \to A \rtimes^{a,(\cup)}_\alpha P$ is an isomorphism. In addition, $A \rta_\alpha P$ and $A \rta_{\alpha,r} P$ can be canonically identified with the ordinary full and reduced crossed product by the group $P$. The reason is that a group does not have any proper (right) ideals, so that both the families $\cJ$ and $\cJ\un$ coincide with the trivial family $\gekl{P,\emptyset}$ in case $P$ is a group.

As we have already mentioned, our construction of semigroup C*-algebras extends the one presented by A. Nica in \cite{Ni}. Let us now explain in detail why this is the case:

A. Nica considers positive cones in so-called quasi-lattice ordered groups. If we reformulate A. Nica's conditions in terms of right ideals, then a quasi-lattice ordered group is a pair $(G,P)$ consisting of a (discrete) subsemigroup $P$ of a (discrete) group $G$ such that $P \cap P^{-1} = \gekl{e}$ where $e$ is the unit element in $G$, and for every $n \geq 1$ and elements $x_1, \dotsc, x_n \in G$,
\bgl
\label{QL}
  P \cap \bigcap_{i=1}^n (x_i \cdot P) \text{ is either empty or of the form } pP \text{ for some } p \in P.
\egl
Note that for $x$ in $G$, we set
\bgl
\label{not-Nica}
  x \cdot P \defeq \menge{xp}{p \in P} \subseteq G.
\egl
Comparing this notation with ours from \eqref{notations}, we obtain that for every $p$, $q$ in $P$, $q^{-1} p P$ in our notation \eqref{notations} is the same as $((q^{-1}p) \cdot P) \cap P$ in notation~\eqref{not-Nica}. More generally (proceeding inductively on $n$), we have for all $p_1, \dotsc, p_n$, $q_1, \dotsc, q_n$ in $P$ that $q_1^{-1} p_1 \dotsm q_n^{-1} p_n P$ in notation~\eqref{notations} coincides with $P \cap (q_1^{-1} p_1) \cdot P \cap \dotsm \cap (q_1^{-1} p_1 \dotsm q_n^{-1} p_n) \cdot P$ in notation~\eqref{not-Nica}. Therefore, for such a semigroup $P$ in a quasi-lattice ordered group $(G,P)$, the family $\cJ$ is simply given by
\bgl
\label{C-Nica}
  \cJ = \menge{pP}{p \in P} \cup \gekl{\emptyset}.
\egl
In other words, the family $\cJ$ consists of the empty set and all principal right ideals of $P$. With this observation, it is now easy to identify A. Nica's construction with ours:

First of all, our definition of the reduced semigroup C*-algebra $C^*_r(P)$ is exactly the same as A. Nica's (see \cite{Ni}, \S~2.4; A. Nica denotes his reduced semigroup C*-algebra by $\cW(G,P)$).

Let us now treat the full versions. A. Nica defines the full semigroup C*-algebra of $P$ (or of the pair $(G,P)$) as the universal C*-algebra for covariant representations of $P$ by isometries. He denotes this C*-algebra by $C^*(G,P)$. To be more precise, this means that $C^*(G,P)$ is the universal C*-algebra generated by isometries $\menge{v(p)}{p \in P}$ subject to the relations
\bglnoz
  && \text{I}_{\text{Nica}} \text{. } v(p) v(q) = v(pq) \\
  && \text{II}_{\text{Nica}} \text{. } v(p) v(p)^* v(q) v(q)^* = 
  \bfa
    v(r) v(r)^* \falls pP \cap qP = rP \text{ for some } r \in P \\
    0 \falls pP \cap qP = \emptyset
  \efa
\eglnoz
for all $p$, $q$ in $P$. Note that by condition~\eqref{QL}, there are only these two possibilities $pP \cap qP = rP$ for some $r \in P$ or $pP \cap qP = \emptyset$.

Now we can construct mutually inverse homomorphisms $C^*(P) \leftrightharpoons C^*(G,P)$ as follows: Send $C^*(P) \ni v_p$ to $v(p) \in C^*(G,P)$ and $C^*(P) \ni e_X$ to $0 \in C^*(G,P)$ if $X = \emptyset$ and to $v(p) v(p)^*$ if $X=pP$ (compare \eqref{C-Nica}). Such a homomorphism $C^*(P) \to C^*(G,P)$ exists as relation I.(i) is exactly relation I$_{\text{Nica}}$ and relation I.(ii) is satisfied as $v_p e_{qP} v_p^* \ma v(p) v(q) v(q)^* v(p)^* \overset{\text{I}_{\text{Nica}}}{=} v(pq) v(pq)^*$ and $e_{pqP} \ma v(pq) v(pq)^*$. Moreover, relations II.(i) and II.(ii) are obviously satisfied, and relation II.(iii) corresponds precisely to relation II$_{\text{Nica}}$. For the homomorphism in the reverse direction, set $C^*(P) \ni v_p \mafr v(p) \in C^*(G,P)$. Such a homomorphism exists because relation I$_{\text{Nica}}$ is relation I.(i), and we have in $C^*(P)$
\bgloz
  v_p v_p^* v_q v_q^* \overset{II.(i)}{=} v_p e_P v_p^* v_q e_P v_q^* \overset{I.(ii)}{=} e_{pP} e_{qP} = \e{pP \cap qP}.
\egloz
If $pP \cap qP$ is of the form $rP$ for some $r$ in $P$, then $e_{pP \cap qP} = e_{rP} = v_r e_P v_r^* = v_r v_r^*$, and if $pP \cap qP = \emptyset$, then $\e{pP \cap qP} = e_{\emptyset} \overset{II.(ii)}{=} 0$. Therefore, relation II$_{\text{Nica}}$ is satisfied. Hence we have seen that $C^*(P)$ and $C^*(G,P)$ are canonically isomorphic. Moreover, we will also see in Corollary~\ref{Nica:C=C-un} that if $P$ is the positive cone in a quasi-lattice ordered group, then the canonical homomorphism $\pi\un: C^*(P) \to {C^*}\un(P)$ is an isomorphism.

So for the special semigroups which A. Nica considers, our constructions indeed coincide with A. Nica's. We refer the reader to \cite{Ni}, Sections~1 and 5 for concrete examples already discussed by A. Nica.

Furthermore, let us compare our construction with the one in \cite{C-D-L}. Given a ring of integers $R$ in a number field, the Toeplitz algebra $\fT[R]$ is defined as the universal C*-algebra generated by
\bglnoz
  && \text{unitaries } \menge{u^b}{b \in R}, \\
  && \text{isometries } \menge{s_a}{a \in R\reg = R \setminus \gekl{0}} \\
  && \text{and projections } \menge{e_I}{(0) \neq I \triangleleft R}
\eglnoz
subject to the relations
\bgln
\label{C-L1}
  &&  u^b s_a u^d s_c = u^{b+ad} s_{ac} \\
\label{C-L2}
  && e_{I \cap J} = e_I \cdot e_J \text{, } e_R = 1 \\
\label{C-L3}
  && s_a e_I s_a^* = e_{aI} \\
\label{C-L4}
  && u^b e_I u^{-b} = e_I \falls b \in I \text{ and } u^b e_I u^{-b} \perp e_I \falls b \notin I.
\egln
Alternatively, we can consider the $ax+b$-semigroup over the ring of integers $R$. It is given by $R \rtimes R\reg = \menge{(b,a)}{b \in R, a \in R\reg}$ where $R\reg = R \setminus \gekl{0}$, and the binary operation is defined by $(b,a)(d,c) = (b+ad,ac)$. Since $R$ is an integral domain, this semigroup $R \rtimes R\reg$ is left cancellative. So we can apply our construction and consider the semigroup C*-algebra $C^*(R \rtimes R\reg)$.

Our goal is to show that $C^*(R \rtimes R\reg)$ and $\fT[R]$ are canonically isomorphic. To see this, we first make two observations:

The relations \eqref{C-L2} and \eqref{C-L4} may be replaced by the stronger relations
\bgln
\label{C-L2i}
  && e_R = 1 \\
\label{C-L2ii}
  && u^b e_I u^{-b} = e_I \fa b \in I \\
\label{C-L2iii}
  && u^{b_1} e_{I_1} u^{-b_1} u^{b_2} e_{I_2} u^{-b_2} = 
  \bfa
    u^d e_{I_1 \cap I_2} u^{-d} \falls (b_1 + I_1) \cap (b_2 + I_2) = d + I_1 \cap I_2 \\
    0 \falls (b_1 + I_1) \cap (b_2 + I_2) = \emptyset.
  \efa
\egln
First of all, it is easy to see that the two cases which appear in \eqref{C-L2iii} are the only possible cases. To see that the relations~\eqref{C-L1}, \eqref{C-L3}, \eqref{C-L2i}--\eqref{C-L2iii} are actually equivalent to the relations~\eqref{C-L1} -- \eqref{C-L4}, we have to prove that the relations \eqref{C-L1} -- \eqref{C-L4} imply \eqref{C-L2iii}. The remaining implications are obvious. Now, if $(b_1+I_1) \cap (b_2+I_2) = \emptyset$, then $-b_1+b_2$ does not lie in $I_1+I_2$. Hence
\bgloz
  u^{b_1} e_{I_1} u^{-b_1} u^{b_2} e_{I_2} u^{-b_2} 
  \overset{\eqref{C-L2}}{=} u^{b_1} e_{I_1} 
  \underbrace{e_{I_1+I_2} u^{-b_1+b_2} e_{I_1+I_2}}_{= \: 0 \text{ by } \eqref{C-L4}}
  e_{I_2} u^{-b_2} = 0.
\egloz
If $(b_1+I_1) \cap (b_2+I_2) = d + I_1 \cap I_2$, then we can find elements $r_1, r_2 \in R$ so that $d=b_1+r_1=b_2+r_2 \Rarr -b_1+b_2 = r_1-r_2$. We conclude that
\bglnoz
  && u^{b_1} e_{I_1} u^{-b_1} u^{b_2} e_{I_2} u^{-b_2} = u^{b_1} e_{I_1} u^{r_1} u^{-r_2} e_{I_2} u^{-b_2} \\
  && \overset{\eqref{C-L4}}{=} u^{b_1} u^{r_1} e_{I_1} e_{I_2} u^{-r_2} u^{b_2} 
  \overset{\eqref{C-L1} \text{, } \eqref{C-L2}}{=} u^d \e{I_1 \cap I_2} u^{-d}.
\eglnoz

Moreover, using the fact that $R$ is a Dedekind domain (the definition of a Dedekind domain is for instance given in \cite{Neu}, Chapter~I, Definition~(3.2)), we can deduce that every ideal $(0) \neq I \triangleleft R$ is of the form $I = ((c^{-1}a) \cdot R) \cap R$ for some $a, c \in R\reg$. (Here $(\cdot)^{-1}$ stands for the inverse in the multiplicative group of the quotient field of $R$.) A proof of this observation is given in \cite{C-D-L}, Lemma~4.15. Here is an alternative proof: Since $R$ is a Dedekind domain, we can find non-zero prime ideals $P_1, \dotsc, P_n$ so that $I = P_1^{\nu_1} \dotsm P_n^{\nu_n}$. By strong approximation (see \cite{Bour2}, Chapitre~VII, \S~2.4, Proposition~2), there are $a, c \in R\reg$ such that
\bgloz
  aR = P_1^{\nu_1} \dotsm P_n^{\nu_n} I_a \text{ for some ideal } I_a \text{ which is coprime to } P_1, \dotsc, P_n
\egloz
and
\bgloz
  cR = I_a I_c \text{ for some ideal } I_c \text{ which is coprime to } I_a \text{ and } P_1, \dotsc, P_n.
\egloz
We then have
\bgloz
  (c^{-1} a) \cdot R = P_1^{\nu_1} \dotsm P_n^{\nu_n} (I_c)^{-1}
\egloz
so that
\bgloz
  ((c^{-1} a) \cdot R) \cap R = P_1^{\nu_1} \dotsm P_n^{\nu_n} = I.
\egloz
This proof shows that in an arbitrary Dedekind domain $R$, every ideal $(0) \neq I \triangleleft R$ is of the form $I = ((c^{-1}a) \cdot R) \cap R$. As $((c^{-1}a) \cdot R) \cap R = c^{-1} (aR)$ where on the right hand side, $c^{-1}$ stands for pre-image (under left multiplication with $c$), it follows that for the semigroup $R \rtimes R\reg$, the family $\cJ$ is given by
\bgloz
  \cJ = \menge{(b+I) \times I\reg}{b \in R, (0) \neq I \triangleleft R} \cup \gekl{\emptyset},
\egloz
where $I\reg = I \cap R\reg = I \setminus \gekl{0}$. Again, this not only holds for rings of integers, but for arbitrary Dedekind domains. 

We can now construct mutually inverse homomorphisms $C^*(R \rtimes R\reg) \leftrightharpoons \fT[R]$ by setting
\bgloz
  v_{(b,a)} \ma u^b s_a \text{, } e_{(b+I) \times I\reg} \ma u^b e_I u^{-b} \text{, } e_{\emptyset} \ma 0 
\egloz
and
\bgloz
  v_{(b,1)} \mafr u^b \text{, } v_{(0,a)} \mafr s_a \text{, } e_{I \times I\reg} \mafr e_I.
\egloz
To see that these homomorphisms really exist, we have to compare the relations from Definition~\ref{full-semigp-C} defining $C^*(R \rtimes R\reg)$ with the relations~\eqref{C-L1}, \eqref{C-L3} and \eqref{C-L2i}--\eqref{C-L2iii}. It is easy to see that
\bglnoz
  && \text{relation I.(i) corresponds to relation } \eqref{C-L1}, \\
  && \text{relation I.(ii) for } p = (0,a) \in R \rtimes R\reg \text{ corresponds to relation } \eqref{C-L3}, \\
  && \text{relation II.(i) is relation } \eqref{C-L2i}, \\
  && \text{relation I.(ii) for } p = (b,1) \in R \rtimes R\reg \text{ is relation } \eqref{C-L2ii} \\
  && \text{and relation II.(iii), together with relation II.(ii), is relation } \eqref{C-L2iii}.
\eglnoz
This proves that $C^*(R \rtimes R\reg)$ and $\fT[R]$ are canonically isomorphic.

\subsection{Functoriality}
At this point, we would like to address the question of functoriality: Given a homomorphism $\varphi: P \to Q$ between left cancellative semigroups, does $\varphi$ induce a homomorphism of the semigroup C*-algebras by the formula $v_p \ma v_{\varphi(p)}$?

It is not clear what the answer to this question in general is because the assignment $v_p \ma v_{\varphi(p)}$ has to be compatible with the extra relations we have built into our constructions. One thing that is clear is that a homomorphism $C^*(P) \to C^*(Q)$ is uniquely determined by the requirement that $v_p$ is sent to $v_{\varphi(p)}$ for all $p$ in $P$. The reason is that $C^*(P)$ is generated as a C*-algebra by the isometries $v_p$ (see Corollary~\ref{v_p-generate}). However, for special semigroups, namely $ax+b$-semigroups over integral domains, we can say more about functoriality. 

We consider the following setting: Let $R$ be an integral domain, i.e. a commutative ring with unit but without zero-divisors. As we did before in the case of rings of integers, we can form the $ax+b$-semigroup $P_R$ over $R$. To be more precise, $P_R$ is the semidirect product $R \rtimes R\reg$, where $R\reg = R \setminus \gekl{0}$ acts multiplicatively on $R$. This means that $P_R = \menge{(b,a)}{b \in R, a \in R\reg}$ and the binary operation is given by $(b,a)(d,c) = (b+ad,ac)$. $P_R$ is left cancellative because $R$ has no zero-divisors. Thus we can form the semigroup C*-algebra $C^*(P_R)$. Let us describe the family $\cJ_{P_R}$ given by \eqref{J-explicit} for this semigroup $P_R$. Given an ideal $I$ of $R$, we denote its image under left multiplication by $a \in R\reg$ by $aI$ and its pre-image under left multiplication with $a \in R\reg$ by $a^{-1} I$, i.e. $aI = \menge{ar}{r \in I}$ and $a^{-1} I = \menge{r \in R}{ar \in I}$. Let $\cI(R)$ be the smallest family of ideals of $R$ which contains $R$, which is closed under left multiplications as well as pre-images under left multiplications, i.e. $a \in R\reg, I \in \cI(R) \Rarr aI, a^{-1} I \in \cI(R)$, and finite intersections, i.e. $I, J \in \cI(R) \Rarr I \cap J \in \cI(R)$. By definition, we have
\bgloz
  \cI(R) = \menge{\bigcap_{j=1}^N (c_{j,1})^{-1} a_{j,1} \dotsm (c_{j,n_j})^{-1} a_{j,n_j} R}{N, n_j \in \Zz \pos; a_{j,k}, c_{j,k} \in R\reg}.
\egloz
We then have
\bgloz
  \cJ_{P_R} = \menge{(b+I) \times I\reg}{b \in R, I \in \cI(R)} \cup \gekl{\emptyset},
\egloz
where $I\reg = I \cap R\reg = I \setminus \gekl{0}$.

Now assume that $S$ is another integral domain, and let $P_S$ be the $ax+b$-semigroup over $S$. Moreover, let $\phi$ be a ring homomorphism $R \to S$. If $\phi$ is injective, it induces a semigroup homomorphism $\varphi: P_R \to P_S$ which sends $P_R \ni (b,a)$ to $(\phi(b),\phi(a)) \in P_S$. Extending the functorial results on Toeplitz algebras associated with rings of integers in number fields from \cite{C-D-L}, Proposition~3.2, we show that there exists a homomorphism $C^*(P_R) \to C^*(P_S)$ sending $v_p$ to $v_{\varphi(p)}$ for every $p \in P$ if $\varphi$ comes from a ring monomorphism $\phi$ such that the quotient $S / \phi(R)$ (in the category of $\phi(R)$-modules) is a flat $\phi(R)$-module.
\blemma
Assume that for all ideals $I$ and $J$ of $R$ which lie in $\cI(R)$, we have
\begin{itemize}
  \item[(a)] $(\phi(I)S) \cap \phi(R) = \phi(I)$
  \item[(b)] $\phi(I)S \cap \phi(J)S = \phi(I \cap J)S$.
\end{itemize}
Then there exists a homomorphism $C^*(P_R) \to C^*(P_S)$ sending $v_p$ to $v_{\varphi(p)}$ for every $p \in P_R$.
\elemma
By $\phi(I)S$, we mean the ideal of $S$ generated by $\phi(I)$.
\bproof
By universal property of $C^*(P_R)$, there exists a homomorphism $C^*(P_R) \to C^*(P_S)$ sending $C^*(P_R) \ni v_p$ to $v_{\varphi(p)} \in C^*(P_S)$ for every $p \in P_R$ and $C^*(P_R) \ni \e{(b+I) \times I\reg}$ to $\e{(\phi(b)+\phi(I)S)\times (\phi(I)S)\reg} \in C^*(P_S)$ for every $b \in R$, $I \in \cI(R)$. To see this, we first of all have to prove that for every $(b+I) \times I\reg \in \cJ_{P_R}$, the right ideal $(\phi(b)+\phi(I)S) \times (\phi(I)S)\reg$ lies in $\cJ_{P_S}$. It suffices to show that for every $I \in \cI(R)$, the ideal $\phi(I)S$ lies in $\cI(S)$, where 
\bgloz
  \cI(S) = \menge{\bigcap_{j=1}^N (c_{j,1})^{-1} a_{j,1} \dotsm (c_{j,n_j})^{-1} a_{j,n_j} S}{N, n_j \in \Zz \pos; a_{j,k}, c_{j,k} \in S\reg}.
\egloz
All we have to prove is that for all $a,c \in R\reg$ and every $I \in \cI(R)$, we have
\bgln
\label{aI}
  \phi(aI)S &=& \phi(a)(\phi(I)S), \\
\label{cI}
  \phi(c^{-1} I)S &=& \phi(c)^{-1} (\phi(I)S).
\egln
\eqref{aI} is obviously true. For \eqref{cI}, we observe that
\bglnoz
  && \phi(c)(\phi(c^{-1} I)S) = \phi(c(c^{-1} I))S = \phi(I \cap cR)S \\
  &\overset{\text{(b)}}{=}& \phi(I)S \cap \phi(cR)S = \phi(I)S \cap \phi(c)S = \phi(c)(\phi(c)^{-1} (\phi(I)S)).
\eglnoz
Applying $\phi(c)^{-1}$ to both sides of this equation yields $\phi(c^{-1} I)S = \phi(c)^{-1} (\phi(I)S)$, as desired.

Moreover, we have to check that the map
\bgloz
  \cJ_{P_R} \ni (b+I) \times I\reg \ma (\phi(b)+\phi(I)S) \times (\phi(I)S)\reg \in \cJ_{P_S}
\egloz
is compatible with left multiplications, taking pre-images under left multiplications and finite intersections. \eqref{aI} and \eqref{cI} imply compatibility with left multiplications and taking pre-images under left multiplications. It remains to prove compatibility with finite intersections. More precisely, we have to show that if
\bgl
\label{1stcase}
  \rukl{(b+I) \times I\reg} \cap \rukl{(d+J) \times J\reg} = \emptyset,
\egl
then
\bgl
\label{1stconcl}
  \rukl{(\phi(b)+\phi(I)S) \times (\phi(I)S)\reg} \cap \rukl{(\phi(d)+\phi(J)S) \times (\phi(J)S)\reg} = \emptyset,
\egl
and if
\bgl
\label{2ndcase}
  \rukl{(b+I) \times I\reg} \cap \rukl{(d+J) \times J\reg} = (r+I \cap J) \times (I \cap J)\reg \text{ for some } r \in R,
\egl
then
\bgln
\label{2ndconcl}
  && \rukl{(\phi(b)+\phi(I)S) \times (\phi(I)S)\reg} \cap \rukl{(\phi(d)+\phi(J)S) \times (\phi(J)S)\reg} \\
  &=& (\phi(r)+\phi(I \cap J)S) \times (\phi(I \cap J)S)\reg. \nonumber
\egln
Now \eqref{1stcase} holds if and only if $(b+I) \cap (d+J) = \emptyset \LRarr b-d \notin I+J$. If the difference $b-d$ does not lie in $I+J$, then $\phi(b)-\phi(d)$ does not lie in 
\bgloz
  \phi(I+J) \overset{\text{(a)}}{=} \phi(I+J)S \cap \phi(R) = (\phi(I)S+\phi(J)S) \cap \phi(R).
\egloz
Hence $\phi(b)-\phi(d)$ does not lie in $\phi(I)S+\phi(J)S$. This implies $(\phi(b)+\phi(I)S) \cap (\phi(d)+\phi(J)S) = \emptyset$, and \eqref{1stconcl} follows. Moreover, \eqref{2ndcase} holds if and only if $(b+I) \cap (d+J) = r+I \cap J \LRarr r \in (b+I) \cap (d+J)$ for some $r \in R$. If $r$ lies in $b+I$, then $\phi(r)$ lies in $\phi(b)+\phi(I)S$. Similarly, $\phi(r)$ lies in $\phi(d)+\phi(J)S$ if $r$ lies in $d+J$. Thus if \eqref{2ndcase} holds, then $\phi(r)$ lies in $(\phi(b)+\phi(I)S) \cap (\phi(d)+\phi(J)S)$. This implies
\bgloz
  (\phi(b)+\phi(I)S) \cap (\phi(d)+\phi(J)S) = \phi(r)+\phi(I)S \cap \phi(J)S \overset{\text{(b)}}{=} \phi(r)+\phi(I \cap J)S.
\egloz
This implies \eqref{2ndconcl}.
\eproof
\bcor
Assume that $\phi: R \to S$ is an inclusion of integral domains such that the quotient $S / \phi(R)$ of the $\phi(R)$-module $S$ by the $\phi(R)$-module $\phi(R)$ (in the category of $\phi(R)$-modules) is a flat $\phi(R)$-module. Let $P_R$ and $P_S$ be the $ax+b$-semigroups over $R$ and $S$, respectively, and let $\varphi: P_R \to P_S$ be the semigroup homomorphism induced by $\phi$. Then there exists a homomorphism $\Phi: C^*(P_R) \to C^*(P_S)$ sending $C^*(P_R) \ni v_p$ to $v_{\varphi(p)} \in C^*(P_S)$.
\ecor
We remark that the condition of flatness already appears in \cite{C-D-L}, Lemma~3.1.
\bproof
If $S / \phi(R)$, the quotient in the category of $\phi(R)$-modules of $S$ by $\phi(R)$, is a flat $\phi(R)$-module, then $S$ itself is a flat $\phi(R)$-module by \cite{Bour1}, Chapitre I, \S~2.5 Proposition~5 using that $\phi(R)$ is flat as a module over itself. Therefore, conditions (a) and (b) from the previous lemma are satisfied, see for instance \cite{Bour1}, Chapitre I, \S~2.6 Proposition~6 and Corollaire (to Proposition~7).
\eproof

\subsection{Comparison of universal C*-algebras}
In the last part of this section, let us compare the universal C*-algebras $C^*(P)$ and ${C^*}\un(P)$. Our goal is to find out under which conditions the canonical homomorphism $\pi\un: C^*(P) \to {C^*}\un(P)$ is an isomorphism. It will be possible to give a criterion in terms of the constructible right ideals of $P$. As a first step, we take a look at the commutative sub-C*-algebras $D(P)$ and $D\un(P)$ of $C^*(P)$ and ${C^*}\un(P)$. Our investigations will also involve the commutative sub-C*-algebra $D_r(P)$ of the reduced semigroup C*-algebra.

\blemma
\label{cond-inj-hom}
Let $D$ be a unital C*-algebra generated by commuting projections $\gekl{f_i}_{i \in I}$. For a non-empty finite set $F \subseteq I$ and a non-empty subset $F' \subseteq F$, define the projection $e(F',F)$ as
\bgloz
  e(F',F) \defeq (\prod_{i \in F'} f_i) \cdot (\prod_{i \in F \setminus F'} (1-f_i)).
\egloz
Then, given a C*-algebra $C$, a homomorphism $\varphi: D \to C$ is injective if and only if for every non-empty finite subset $F \subseteq I$ and $\emptyset \neq F' \subseteq F$ as above,
\bgl
\label{phi-inj}
  \varphi(e(F',F)) = 0 \text{ in } C \text{ implies } e(F',F) = 0 \text{ in } D.
\egl
\elemma
\bproof
If $\varphi$ is injective, then certainly $\varphi(e(F',F)) = 0$ must imply $e(F',F) = 0$. To prove the reverse implication, we set $D_F \defeq C^*(\menge{f_i}{i \in F}) \subseteq D$ for every non-empty finite subset $F \subseteq I$. The non-empty finite subsets of $I$ are ordered by inclusion, and we obviously have
\bgloz
  D = \overline{\bigcup_{\emptyset \neq F \subseteq I \text{ finite}} D_F}.
\egloz
So it remains to prove that if condition~\eqref{phi-inj} holds for a non-empty finite subset $F \subseteq I$, then $\varphi \vert_{D_F}$ is injective.

But since the projections $\menge{f_i}{i \in F}$ commute, it is clear that the projections $e(F',F)$, $\emptyset \neq F' \subseteq F$ are pairwise orthogonal. This implies that
\bgloz
  D_F = \bigoplus_{\emptyset \neq F' \subseteq F} \Cz \cdot e(F',F).
\egloz
Hence it follows that $\varphi \vert_{D_F}$ is injective if and only if \eqref{phi-inj} holds for every non-empty subset $F'$ of $F$.
\eproof
As a next step, we work out how the projections $e(F',F)$ look like in the following situation: Let $D = D\un(P)$, $I = \cJ\un$ and for every $X \in \cJ\un$, set $f_X \defeq e_X \in {C^*}\un(P)$ (see Definition~\ref{full-semigp-C-union}).
\blemma
\label{e=e_X-e_Y}
For every non-empty finite subset $F \subseteq \cJ\un$ and every $\emptyset \neq F' \subseteq F$, there exist $X, Y \in \cJ\un$ with $Y \subseteq X$ such that $e(F',F) = e_X - e_Y$.
\elemma
\bproof
Let us proceed inductively on $\abs{F}$. The starting point $\abs{F}=1$ is trivial. We assume that the claim is proven whenever $\abs{F}=n$. Let $F$ be a finite subset of $\cJ\un$ with $\abs{F}=n+1$. If $F' = F$ then our assertion obviously follows from relation II$\un$.(iii). If $\emptyset \neq F' \subsetneq F$, then we can find a subset $F_n$ of $\cJ\un$ with $\abs{F_n}=n$ and $F' \subseteq F_n \subseteq F$. Let $F = F_n \cup \gekl{X_{n+1}}$. We know by induction hypothesis that there exist $X_n, Y_n \in \cJ\un$ with $Y_n \subseteq X_n$ such that $e(F',F_n) = e_{X_n} - e_{Y_n}$. Therefore,
\bglnoz
  && e(F',F) = e(F',F_n)(1-e_{X_{n+1}}) = (e_{X_n}-e_{Y_n})(1-e_{X_{n+1}}) \\
  &\overset{\text{II}\un.(iii)}{=}& e_{X_n}-e_{Y_n}-\e{X_n \cap X_{n+1}}+\e{Y_n \cap X_{n+1}} 
  \overset{\text{II}\un.(iv)}{=} e_{X_n}-\e{Y_n \cup (X_n \cap X_{n+1})}.
\eglnoz
Set $X=X_n$, $Y=Y_n \cup (X_n \cap X_{n+1})$ and we are done.
\eproof
\bcor
\label{lambda/D-union}
$\lambda\un \vert_{D\un(P)}: D\un(P) \to D_r(P)$ is an isomorphism.
\ecor
\bproof
It is clear that $\lambda\un \vert_{D\un(P)}$ is surjective, thus it remains to prove injectivity. We want to apply Lemma~\ref{cond-inj-hom} to $D = D\un(P) = C^*(\menge{e_X}{X \in \cJ\un})$, $C = D_r(P)$ and $\varphi = \lambda\un \vert_{D\un(P)}$. For a non-empty finite subset $F \subseteq \cJ\un$ and $\emptyset \neq F' \subseteq F$, Lemma~\ref{e=e_X-e_Y} tells us that there are $X, Y \in \cJ\un$ with $Y \subseteq X$ such that $e(F',F) = e_X - e_Y$. Now $\lambda\un(e_X - e_Y) = E_X - E_Y$, and $E_X - E_Y$ vanishes as an operator on $\ell^2(P)$ if and only if $X=Y$. But $X=Y$ obviously implies $e(F',F) = e_X - e_Y = 0$ in $D\un(P)$. Therefore, Lemma~\ref{cond-inj-hom} implies that $\lambda\un \vert_{D\un(P)}$ must be injective.
\eproof
\bcor
\label{C-un(product)2}
Given two left cancellative semigroups $P$ and $Q$, we can identify ${C^*}\un(P \times Q)$ with ${C^*}\un(P) \otimes_{\max} {C^*}\un(Q)$ via a homomorphism sending $v_{(p,q)}$ to $v_p \otimes v_q$ for every $p \in P$ and $q \in Q$.
\ecor
\bproof
As explained in Remark~\ref{C-un(product)1}, all we have to do is to construct a homomorphism $D\un(P \times Q) \to {C^*}\un(P) \otimes_{\max} {C^*}\un(Q)$ which sends for all $X \in \cJ_P$ and $Y \in \cJ_Q$ the projection $e_{X \times Y}$ to $e_X \otimes e_Y$. But we know by the previous lemma that $D\un(P \times Q) \cong D_r(P \times Q)$, $D\un(P) \cong D_r(P)$ and $D\un(Q) \cong D_r(Q)$. Moreover, the isomorphism $C^*_r(P \times Q) \cong C^*_r(P) \otimes_{\min} C^*_r(Q)$ from Lemma~\ref{C(product)} obviously identifies $D_r(P \times Q)$ with $D_r(P) \otimes_{\min} D_r(Q)$. Thus the desired homomorphism is given by
\bglnoz
  && D\un(P \times Q) \cong D_r(P \times Q) \cong D_r(P) \otimes_{\min} D_r(Q) \cong D_r(P) \otimes_{\max} D_r(Q) \\
  &\cong& D\un(P) \otimes_{\max} D\un(Q) \to {C^*}\un(P) \otimes_{\max} {C^*}\un(Q).
\eglnoz
\eproof
Now we come to the main result concluding this circle of ideas.
\bprop
\label{pi-iso}
The following statements are equivalent:
\begin{itemize}
\item[(i)] If $X = \bigcup_{j=1}^n X_j$ for $X, X_1, \dotsc, X_n \in \cJ$, then $X = X_j$ for some $1 \leq j \leq n$.
\item[(ii)] $\pi\un \vert_{D(P)}: D(P) \to D\un(P)$ is an isomorphism.
\item[(iii)] $\pi\un: C^*(P) \to {C^*}\un(P)$ is an isomorphism.
\item[(iv)] There exists a homomorphism $\Delta\un: {C^*}\un(P) \to {C^*}\un(P) \otimes_{\max} {C^*}\un(P)$ which sends (for all $p \in P$) $v_p$ to $v_p \otimes v_p$.
\item[(v)] There exists a homomorphism $\Delta\un_D: D\un(P) \to D\un(P) \otimes_{\max} D\un(P)$ which sends (for all $X \in \cJ$) $e_X$ to $e_X \otimes e_X$.
\end{itemize}
\eprop
\bproof
\an{(i) $\Rarr$ (ii)}: Since by Corollary~\ref{lambda/D-union}, $\lambda\un \vert_{D\un(P)}$ is an isomorphism and because we always have $\lambda = \lambda\un \circ \pi\un$, statement (ii) is equivalent to \an{$\lambda \vert_{D(P)}$ is an isomorphism}. $\lambda \vert_{D(P)}$ is obviously surjective, so it remains to prove injectivity. We want to apply Lemma~\ref{cond-inj-hom} to $D=D(P)$, $I=\cJ$, $f_X \defeq e_X \in D(P)$ for $X \in \cJ$, $C = D_r(P)$ and $\varphi = \lambda \vert_{D(P)}$. Given a non-empty finite subset $F \subseteq \cJ$ and $\emptyset \neq F' \subseteq F$, it is immediate that $\lambda(e(F',F)) = \E{(\bigcap_{X' \in F'} X') \setminus (\bigcup_{Y \in F \setminus F'} Y)}$ where $\E{(\bigcap_{X' \in F'} X') \setminus (\bigcup_{Y \in F \setminus F'} Y)}$ is the orthogonal projection onto the subspace
\bgloz
  \ell^2 \rukl{(\bigcap_{X' \in F'} X') \setminus (\bigcup_{Y \in F \setminus F'} Y)} \subseteq \ell^2(P).
\egloz
Assume that $\lambda(e(F',F))$ vanishes. Then $X \defeq \bigcap_{X' \in F'} X'$ must be a subset of $\bigcup_{Y \in F \setminus F'} Y$. Now $X$ lies in $\cJ$, and $X \subseteq \bigcup_{Y \in F \setminus F'} Y$ implies $X = \bigcup_{Y \in F \setminus F'} (Y \cap X)$. But statement (i) tells us that this can only happen if there exists $Y \in F \setminus F'$ with $Y \cap X = X$, or equivalently, $X \subseteq Y$. Thus $e_X = e_{X \cap Y} \overset{II.(iii)}{=} e_X \cdot e_Y$, and we conclude that $e_X(1-e_Y)=0$. Hence it follows that
\bgloz
  e(F',F) = e_X(1-e_Y) \cdot \prod_{Y \neq Z \in F \setminus F'} (1-e_Z) = 0.
\egloz
So we have seen that condition~\eqref{phi-inj} holds. Therefore $\lambda \vert_{D(P)}$ is injective.

\an{(ii) $\Rarr$ (iii)}: This follows from the crossed product descriptions of $C^*(P)$ and ${C^*}\un(P)$ from Lemma~\ref{cropro-descriptions} and the fact that $\pi\un \vert_{D(P)}$ is $P$-equivariant with respect to the actions $\tau$ and $\tau\un$.

\an{(iii) $\Rarr$ (iv)}: It follows from universal property of $C^*(P)$ that there exists a homomorphism $\Delta: C^*(P) \to C^*(P) \otimes_{\max} C^*(P)$ which sends $v_p$ to $v_p \otimes v_p \in C^*(P) \odot C^*(P) \subseteq C^*(P) \otimes_{\max} C^*(P)$ and $e_X$ to $e_X \otimes e_X \in C^*(P) \odot C^*(P) \subseteq C^*(P) \otimes_{\max} C^*(P)$ for every $p \in P$ and $X \in \cJ$. The reason is that relations I and II are obviously valid with $v_p \otimes v_p$ in place of $v_p$ and $e_X \otimes e_X$ in place of $e_X$. Now set $\Delta\un \defeq ((\pi\un)^{-1} \otimes_{\max} (\pi\un)^{-1}) \circ \Delta \circ \pi\un$.

\an{(iv) $\Rarr$ (v)}: Just restrict $\Delta\un$ to $D\un(P)$, i.e. set $\Delta\un_D \defeq \Delta\un \vert_{D\un(P)}$.

\an{(v) $\Rarr$ (i)}:
Let $D$ be the sub-*-algebra of $D\un(P)$ generated by the projections $\menge{e_X}{X \in \cJ\un}$. By relation II$\un$.(iii), the set $\menge{e_X}{X \in \cJ}$ is multiplicatively closed, and by relation II$\un$.(iv), $D = \lspan(\menge{e_X}{\emptyset \neq X \in \cJ})$. Restricting $\Delta\un_D$ to $D$, we obtain a homomorphism $D \to D \odot D$ which is determined by $e_X \ma e_X \otimes e_X$ for every $X \in \cJ$. Let us denote this restriction by $\Delta_D$.

We can now deduce from the existence of such a homomorphism $\Delta_D$ that the set $\menge{e_X}{\emptyset \neq X \in \cJ}$ is a $\Cz$-basis of $D$. As $\menge{e_X}{\emptyset \neq X \in \cJ}$ generates $D$ as a $\Cz$-vector space, we can always find a subset $\ti{\cJ}$ of $\cJ \setminus \gekl{\emptyset}$ such that $\menge{e_X}{X \in \ti{\cJ}}$ is a $\Cz$-basis for $D$. It then follows that $\menge{e_{\ti{X}} \otimes e_{\ti{Y}}}{\ti{X}, \ti{Y} \in \ti{\cJ}}$ is a $\Cz$-basis of $D \odot D$.

Now take $\emptyset \neq X \in \cJ$. We can find finite subsets $\gekl{X_i} \subseteq \ti{\cJ}$ and $\gekl{\alpha_i} \subseteq \Cz$ with $e_X = \sum_i \alpha_i e_{X_i}$. Applying $\Delta_D$ yields
\bgloz
  \sum_{i,j} \alpha_i \alpha_j e_{X_i} \otimes e_{X_j} = e_X \otimes e_X = \Delta_D(e_X) 
  = \sum_i \alpha_i \Delta_D(e_{X_i}) = \sum_i \alpha_i e_{X_i} \otimes e_{X_i}. 
\egloz
Hence it follows that among the $\alpha_i$s, there can only be one non-zero coefficient which must be $1$. The corresponding vector $e_{X_i}$ must then coincide with $e_X$. This implies $e_X \in \menge{e_{\ti{X}}}{\ti{X} \in \ti{\cJ}}$, i.e. $\menge{e_X}{\emptyset \neq X \in \cJ} = \menge{e_{\ti{X}}}{\ti{X} \in \ti{\cJ}}$ is a $\Cz$-basis of $D$.

Now assume that there are $X, X_1, \dotsc, X_n \in \cJ$ with $X = \bigcup_{j=1}^n X_j$. We necessarily have $X_j \subseteq X$ for all $1 \leq j \leq n$. Moreover, $X_j \subsetneq X$ implies $e_{X_j} \lneq e_X$ because $\lambda(e_{X_j}) = E_{X_j} \lneq E_X = \lambda(e_X)$ as concrete operators on $\ell^2(P)$. Using relation II$\un$.(iv), we obtain from $X = \bigcup_{j=1}^n X_j$ that
\bgl
\label{rel:X-X_j}
  e_X = \sum_{\emptyset \neq F \subseteq \gekl{1, \dotsc, n}} (-1)^{\abs{F}+1} \e{\bigcap_{j \in F} X_j}
\egl
holds in $D$. But if all the $X_j$s ($1 \leq j \leq n$) are strictly contained in $X$, then \eqref{rel:X-X_j} would give a non-trivial relation among $e_X$ and those projections $\e{\bigcap_{j \in F} X_j}$, $\emptyset \neq F \subseteq \gekl{1, \dotsc, n}$ which are non-zero. But this contradicts our observation that $\menge{e_X}{\emptyset \neq X \in \cJ}$ is a $\Cz$-basis of $D$. Hence we conclude that one of the $X_j$s must be equal to $X$. This proves (i).
\eproof
\bremark
This proposition does not really have much to do with semigroups. It actually is a statement about families of subsets of a fixed set and a projection-valued spectral measure defined on this family.
\eremark
\bdefin
\label{ind}
We call $\cJ$ independent (or we also say that the constructible right ideals of $P$ are independent) if the right ideals in $\cJ$ satisfy (i) from Proposition~\ref{pi-iso}.
\edefin
Note that statement (i) is equivalent to the following one: For all $X$, $X_1$, ..., $X_n$ in $\cJ$ such that $X_1$, ..., $X_n$ are proper subsets of $X$ ($X_i \subsetneq X$ for all $1 \leq i \leq n$), then $\bigcup_{i=1}^n X_i$ must be a proper subset of $X$ ($\bigcup_{i=1}^n X_i \subsetneq X$).

\bcor
\label{un-cond->lrr/D}
The constructible right ideals of $P$ are independent if and only if the restriction of the left regular representation to the commutative sub-C*-algebra $D(P)$ of the full semigroup C*-algebra $C^*(P)$ is an isomorphism.
\ecor
\bproof
This follows immediately from the equivalence of (i) and (ii) in Proposition~\ref{pi-iso} and from Corollary~\ref{lambda/D-union}.
\eproof
An immediate question that comes to mind after Proposition~\ref{pi-iso} is which semigroups have independent constructible right ideals. The general answer is not known to the author. But we can discuss two particular cases:
\blemma
The constructible right ideals of the positive cone in a quasi-lattice ordered group are independent.
\elemma
\bproof
This follows immediately from the observation that for a semigroup $P$ which is the positive cone in a quasi-lattice ordered group, the family $\cJ$ consists of the empty set and all principal right ideals of $P$, see \eqref{C-Nica}.
\eproof
As an immediate consequence of this lemma and Proposition~\ref{pi-iso}, we obtain
\bcor
\label{Nica:C=C-un}
If $P$ is the positive cone in a quasi-lattice ordered group, then the canonical homomorphism $\pi\un: C^*(P) \to {C^*}\un(P)$ is an isomorphism.
\ecor
Another class of semigroups with independent constructible right ideals is given as follows:
\blemma
\label{Dedekind->union}
Let $R$ be a Dedekind domain. Then the constructible right ideals of the $ax+b$-semigroup $P_R$ over $R$ are independent .
\elemma
\bproof
Recall that we have shown above when we identified Toeplitz algebras of rings of integers with full semigroup C*-algebras of the corresponding $ax+b$-semigroups that
\bgloz
  \cJ_{P_R} = \menge{(b+I) \times I\reg}{b \in R, (0) \neq I \triangleleft R} \cup \gekl{\emptyset}.
\egloz
Assume that we have
\bgloz
  (b+I) \times I\reg = \bigcup_{j=1}^n (b_j+I_j) \times I_j\reg
\egloz
with $(b_j+I_j) \times I_j\reg \subsetneq (b+I) \times I\reg$ for all $1 \leq j \leq n$. Then it follows that $I = \bigcup_{j=1}^n I_j$ with $I_j \subsetneq I$ for all $1 \leq j \leq n$.

Because $R$ is a Dedekind domain, we can find non-zero prime ideals $P_1$, ..., $P_N$ of $R$ so that
\bgloz
  I = P_1^{\nu_1} \dotsm P_M^{\nu_M} \text{ for some } M \leq N \text{ and } \nu_1, \dotsc, \nu_M > 0
\egloz
and
\bgloz
  I_j = P_1^{\nu_{1,j}} \dotsm P_M^{\nu_{M,j}} \dotsm P_N^{\nu_{N,j}} \text{ for some } \nu_{i,j} \geq 0 \text{ with } \nu_{i,j} \geq \nu_i \fa 1 \leq i \leq M.
\egloz
By strong approximation (see \cite{Bour2}, Chapitre~VII, \S~2.4, Proposition~2), there exists $x \in R$ with the properties
\begin{itemize}
\item[(*)] $x \in P_i^{\nu_i} \setminus P_i^{\nu_i+1}$ for all $1 \leq i \leq M$
\item[(**)] $x \notin P_i$ for all $M < i \leq N$.
\end{itemize}
(*) implies that $x$ lies in $I$. But $x$ does not lie in $I_j$ for any $1 \leq j \leq n$: If $I_j \subseteq P_i$ for some $M < i \leq N$, then (**) implies that $x \notin I_j \subseteq P_i$. If $I_j$ is coprime to $P_i$ for all $M < i \leq N$ (i.e. $\nu_{i,j} = 0$ for all $M < i \leq N$), then $I_j \subsetneq I$ implies $\nu_{i,j} > \nu_i$ for some $1 \leq i \leq M$. So (*) implies that $x \notin I_j \subseteq P_i^{\nu_{i,j}} \subseteq P_i^{\nu_i+1}$. But this implies that $I \subsetneq \bigcup_{j=1}^n I_j$ which contradicts our assumption.
\eproof
In particular, the constructible right ideals of the $ax+b$-semigroup $P_R$ over the ring of integers $R$ in a number field are independent. So by Corollary~\ref{un-cond->lrr/D}, the left regular representation restricted to the commutative sub-C*-algebra $D(P_R)$ is an isomorphism. This explains Corollary~4.16 in \cite{C-D-L} ($\fT[R]$ in \cite{C-D-L} is canonically isomorphic to $C^*(P_R)$ as explained above, and $\fT$ in \cite{C-D-L} is $C^*_r(P_R)$).

\bremark
In the proof of Lemma~\ref{Dedekind->union}, we have just shown that whenever given non-zero ideals $I$, $I_1$, ..., $I_n$ of a Dedekind domain $R$ such that $I_1$, ..., $I_n$ are proper subsets of $I$, then $\bigcup_{i=1}^n I_i$ is a proper subset of $I$. This means that already the non-zero ideals of a Dedekind domain are independent.
\eremark

\section{A variant of our construction for subsemigroups of groups}
\label{variant}

Given a subsemigroup of a group, let us now modify our construction of full semigroup C*-algebras. We impose extra relations besides the ones from Definition~\ref{full-semigp-C}. These relations are motivated by the following
\blemma
\label{iii_G-red}
Let $P$ be a subsemigroup of a group $G$. Given $p_1$, $q_1$, ..., $p_m$, $q_m$ in $P$ with $p_1^{-1} q_1 \dotsm p_m^{-1} q_m = e$ in $G$, then $V_{p_1}^* V_{q_1} \dotsm V_{p_m}^* V_{q_m} = \E{q_m^{-1} p_m \dotsm q_1^{-1} p_1 P}$ in $C^*_r(P)$.
\elemma
\bproof
For $x \in P$, we have $\E{q_m^{-1} p_m \dotsm q_1^{-1} p_1 P} \ve_x = \ve_x$ if $x \in p_m^{-1} p_m \dotsm q_1^{-1} p_1 P$ and $\E{q_m^{-1} p_m \dotsm q_1^{-1} p_1 P} \ve_x = 0$ if $x \notin q_m^{-1} p_m \dotsm q_1^{-1} p_1 P$. A direct computation yields that $(V_{p_1}^* V_{q_1} \dotsm V_{p_m}^* V_{q_m})(\ve_x) \neq 0$ if and only if $x$ lies in $q_m^{-1} p_m \dotsm q_1^{-1} p_1 P$, and in this case, we have $(V_{p_1}^* V_{q_1} \dotsm V_{p_m}^* V_{q_m})(\ve_x) = \ve_{p_1^{-1} q_1 \dotsm p_m^{-1} q_m x} = \ve_x$.
\eproof

\bdefin
\label{full-semigp-C_s}
Let $P$ be a subsemigroup of a group $G$. We let $C^*_s(P)$ be the universal C*-algebra generated by isometries $\menge{v_p}{p \in P}$ and projections $\menge{e_X}{X \in \cJ}$ satisfying the following relations:
\begin{itemize}
\item[I.] $v_{pq} = v_p v_q$,
\item[II.] $e_{\emptyset} = 0$,
\item[III$_G$.] whenever $p_1, q_1, \dotsc, p_m, q_m \in P$ satisfy $p_1^{-1} q_1 \dotsm p_m^{-1} q_m = e$ in $G$, then $$v_{p_1}^* v_{q_1} \dotsm v_{p_m}^* v_{q_m} = \e{q_m^{-1} p_m \dotsm q_1^{-1} p_1 P}$$
\end{itemize}
for all $p$, $q$ in $P$ and $X$, $Y$ in $\cJ$.

As before, we set $D_s(P) \defeq C^*(\menge{e_X}{X \in \cJ}) \subseteq C^*_s(P)$.
\edefin

By universal property of $C^*_s(P)$ and Lemma~\ref{iii_G-red}, there exists a homomorphism $\lambda: C^*_s(P) \to C^*_r(P)$ determined by $\lambda(v_p) = V_p$ and $\lambda(e_X) = E_X$. In particular, $C^*_s(P)$ is non-zero.

It turns out that relation III$_G$ implies the relations I.(ii), II.(i) and II.(iii) from Definition~\ref{full-semigp-C}. Here is an equivalent way of formulating this:
\blemma
\label{C-C_s}
There is a surjective homomorphism $\pi_s: C^*(P) \to C^*_s(P)$ sending $C^*(P) \ni v_p$ to $v_p \in C^*_s(P)$ and $C^*(P) \ni e_X$ to $e_X \in C^*_s(P)$.
\elemma
\bproof
It suffices to check that such a homomorphism exists. We have to show that the relations I.(ii), II.(i) and II.(iii) from Definition~\ref{full-semigp-C} are satisfied in $C^*_s(P)$. The universal property of $C^*(P)$ will then imply existence of $\pi_s$.

II.(i) holds in $C^*_s(P)$ as $e_P \overset{\text{III}_G}{=} v_e^* v_e = 1$. To proceed, we first prove a general result about the family of constructible right ideals of $P$, namely, that it is automatic that $\cJ$ is closed under finite intersections, i.e.
\bgl
\label{J-without-int}
  \cJ = \menge{q_1^{-1} p_1 \dotsm q_m^{-1} p_m P}{m \geq 1; p_i, q_i \in P} \cup \gekl{\emptyset}.
\egl
To prove \eqref{J-without-int}, we first show that for every $p_i, q_i \in P$ and every subset $X$ of $P$,
\bgl
\label{qpint}
  q_1^{-1} p_1 \dotsm q_m^{-1} p_m p_m^{-1} q_m \dotsm p_1^{-1} q_1 X = (q_1^{-1} p_1 \dotsm q_m^{-1} p_m P) \cap X.
\egl
We proceed inductively on $m$:

\an{$m=1$}:
\bgl
\label{m=1}
  q_1^{-1} p_1 p_1^{-1} q_1 X = q_1^{-1} ((p_1 P) \cap q_1 X) = (q_1^{-1} p_1 P) \cap X.
\egl
\an{$m \to m+1$}:
\bglnoz
  && q_1^{-1} p_1 \dotsm q_{m+1}^{-1} p_{m+1} p_{m+1}^{-1} q_{m+1} \dotsm p_1^{-1} q_1 X \\
  &=& (q_1^{-1} p_1 \dotsm q_m^{-1} p_m) (q_{m+1}^{-1} p_{m+1} p_{m+1}^{-1} q_{m+1} (p_m^{-1} q_m \dotsm p_1^{-1} q_1 X)) \\
  &\overset{\eqref{m=1}}{=}& (q_1^{-1} p_1 \dotsm q_m^{-1} p_m) ((q_{m+1}^{-1} p_{m+1} P) \cap (p_m^{-1} q_m \dotsm p_1^{-1} q_1 X)) \\
  &=& (q_1^{-1} p_1 \dotsm q_{m+1}^{-1} p_{m+1} P) \cap (q_1^{-1} p_1 \dotsm q_m^{-1} p_m p_m^{-1} q_m \dotsm p_1^{-1} q_1 X) \\
  &=& (q_1^{-1} p_1 \dotsm q_{m+1}^{-1} p_{m+1} P) \cap (q_1^{-1} p_1 \dotsm q_m^{-1} p_m P) \cap X \text{ (by induction hypothesis)} \\
  &=& (q_1^{-1} p_1 \dotsm q_{m+1}^{-1} p_{m+1} P) \cap X \text{ (as } q_1^{-1} p_1 \dotsm q_{m+1}^{-1} p_{m+1} P \subseteq q_1^{-1} p_1 \dotsm q_m^{-1} p_m P \text{)}.
\eglnoz
This proves \eqref{qpint}.

We deduce that the right hand side in \eqref{J-without-int} is closed under finite intersections. This implies by definition of $\cJ$ that \an{$\subseteq$} in \eqref{J-without-int} holds. As \an{$\supseteq$} obviously holds as well, we have proven \eqref{J-without-int}.

Let us now show that I.(ii) and II.(iii) from Definition~\ref{full-semigp-C} are satisfied in $C^*_s(P)$. As a special case of \eqref{qpint} ($X=P$), we obtain
\bgl
\label{qpint'}
  q_1^{-1} p_1 \dotsm q_m^{-1} p_m p_m^{-1} q_m \dotsm p_1^{-1} q_1 P = q_1^{-1} p_1 \dotsm q_m^{-1} p_m P.
\egl 
Take $p \in P$, $X = q_1^{-1} p_1 \dotsm q_m^{-1} p_m P \in \cJ$ and $Y = s_1^{-1} r_1 \dotsm s_n^{-1} r_n P \in \cJ$. Then
\bglnoz
  && v_p e_X v_p^* = v_p \e{q_1^{-1} p_1 \dotsm q_m^{-1} p_m P} v_p^* 
  \overset{\eqref{qpint'}}{=} v_p \e{q_1^{-1} p_1 \dotsm q_m^{-1} p_m p_m^{-1} q_m \dotsm p_1^{-1} q_1 P} v_p^* \\
  &\overset{\text{III}_G}{=}& v_p v_{q_1}^* v_{p_1} \dotsm v_{q_m}^* v_{p_m} v_{p_m}^* v_{q_m} \dotsm v_{p_1}^* v_{q_1} v_p^* 
  \overset{\text{III}_G}{=} \e{p q_1^{-1} p_1 \dotsm q_m^{-1} p_m p_m^{-1} q_m \dotsm p_1^{-1} q_1 p^{-1} P} \\
  &\overset{\eqref{qpint'}}{=}& \e{p q_1^{-1} p_1 \dotsm q_m^{-1} p_m P} = e_{pX}.
\eglnoz
This proves I.(ii). Moreover,
\bglnoz
  && e_X e_Y = \e{q_1^{-1} p_1 \dotsm q_m^{-1} p_m P} \e{s_1^{-1} r_1 \dotsm s_n^{-1} r_n P} \\
  &\overset{\eqref{qpint'}}{=}& \e{q_1^{-1} p_1 \dotsm q_m^{-1} p_m p_m^{-1} q_m \dotsm p_1^{-1} q_1 P} 
  \e{s_1^{-1} r_1 \dotsm s_n^{-1} r_n r_n^{-1} s_n \dotsm r_1^{-1} s_1 P} \\
  &\overset{\text{III}_G}{=}& 
  v_{q_1}^* v_{p_1} \dotsm v_{q_m}^* v_{p_m} v_{p_m}^* v_{q_m} \dotsm v_{p_1}^* v_{q_1} 
  v_{s_1}^* v_{r_1} \dotsm v_{s_n}^* v_{r_n} v_{r_n}^* v_{s_n} \dotsm v_{r_1}^* v_{s_1} \\
  &\overset{\text{III}_G}{=}&
  \e{s_1^{-1} r_1 \dotsm s_n^{-1} r_n r_n^{-1} s_n \dotsm r_1^{-1} s_1 (q_1^{-1} p_1 \dotsm q_m^{-1} p_m p_m^{-1} q_m \dotsm p_1^{-1} q_1 P)} \\
  &\overset{\eqref{qpint}}{=}& \e{(s_1^{-1} r_1 \dotsm s_n^{-1} r_n P) \cap (q_1^{-1} p_1 \dotsm q_m^{-1} p_m p_m^{-1} q_m \dotsm p_1^{-1} q_1 P)} \\
  &\overset{\eqref{qpint'}}{=}& \e{(s_1^{-1} r_1 \dotsm s_n^{-1} r_n P) \cap (q_1^{-1} p_1 \dotsm q_m^{-1} p_m P)} = e_{X \cap Y}.
\eglnoz
Thus II.(iii) also holds in $C^*_s(P)$.
\eproof
It follows from Corollary~\ref{v_p-generate} that $C^*_s(P)$ is generated by the isometries $\menge{v_p}{p \in P}$. By construction, we have a commutative triangle
\bgloz
\xymatrix{
C^*(P) \ar[d]_{\pi_s} \ar[dr]^{\lambda}
&
\\
C^*_s(P) \ar[r]_{\lambda} & C^*_r(P).
}
\egloz
Since $\pi_s(D(P)) = D_s(P)$, we can restrict this triangle to $D(P)$ and obtain another commutative diagram
\bgloz
\xymatrix{
D(P) \ar[d]_{\pi_s} \ar[dr]^{\lambda}
&
\\
D_s(P) \ar[r]_{\lambda} & D_r(P).
}
\egloz
As $\pi_s: D(P) \to D_s(P)$ is surjective, we deduce from Corollary~\ref{un-cond->lrr/D}
\bcor
If the constructible right ideals of $P$ are independent, then $\lambda \vert_{D_s(P)}: D_s(P) \to D_r(P)$ is an isomorphism.
\ecor
Moreover, we obtain by universal property of $C^*_s(P)$ a homomorphism
\bgl
\label{Delta-C_s}
  \Delta: C^*_s(P) \to C^*_s(P) \otimes_{\max} C^*_s(P), \: v_p \ma v_p \otimes v_p, \: e_X \ma e_X \otimes e_X.
\egl

In the definition of $C^*_s(P)$, we have used the inclusion $P \subseteq G$. However, the C*-algebra $C^*_s(P)$ is independent from $G$ (up to canonical isomorphism). Namely, $C^*_s(P)$ can be viewed as $C^*(P)$ with the extra relations III$_G$ by Lemma~\ref{C-C_s}. To show independence, let $P \subseteq G_1$ and $P \subseteq G_2$ be two embeddings. We want to see that III$_{G_1}$ and III$_{G_2}$ give the same relations. As we do not add relations if $\e{q_m^{-1} p_m \dotsm q_1^{-1} p_1 P} = 0$ in III$_G$, all we have to show is that for all $p_1$, $q_1$, ..., $p_m$, $q_m$ in $P$,
\bgln
\label{C_s-ind}
  && p_1^{-1} q_1 \dotsm p_m^{-1} q_m = e \text{ in } G_1 \text{ and } \e{q_m^{-1} p_m \dotsm q_1^{-1} p_1 P} \neq 0 \\
  &\LRarr& p_1^{-1} q_1 \dotsm p_m^{-1} q_m = e \text{ in } G_2 \text{ and } \e{q_m^{-1} p_m \dotsm q_1^{-1} p_1 P} \neq 0 .\nonumber
\egln
Once this is proven, we conclude that $C^*_s(P)$ is independent from the group into which we embed $P$. By symmetry, it suffices to prove \an{$\Rarr$}. Take $p_1$, $q_1$, ..., $p_m$, $q_m$ in $P$ such that $p_1^{-1} q_1 \dotsm p_m^{-1} q_m = e$ in $G_1$ and $\e{q_m^{-1} p_m \dotsm q_1^{-1} p_1 P} \neq 0$. As the latter condition implies $q_m^{-1} p_m \dotsm q_1^{-1} p_1 P \neq \emptyset$, we can choose $x \in q_m^{-1} p_m \dotsm q_1^{-1} p_1 P$. Then on $\ell^2(P)$, we have $(V_{p_1}^* V_{q_1} \dotsm V_{p_m}^* V_{q_m})(\ve_x) = \ve_x$ as $p_1^{-1} q_1 \dotsm p_m^{-1} q_m = e$ in $G_1$. But we also have $(V_{p_1}^* V_{q_1} \dotsm V_{p_m}^* V_{q_m})(\ve_x) = \ve_{p_1^{-1} q_1 \dotsm p_m^{-1} q_m x}$ where this time, the product $p_1^{-1} q_1 \dotsm p_m^{-1} q_m x$ is taken in $G_2$. Thus we have $p_1^{-1} q_1 \dotsm p_m^{-1} q_m x = x$ in $G_2$, hence $p_1^{-1} q_1 \dotsm p_m^{-1} q_m = e$ in $G_2$. This proves \eqref{C_s-ind}.

We remark that we can also define ${C^*_s}\un(P)$ (see \S~\ref{constructions}) and crossed products $A \rta_{\alpha,s}P$ as in \S~\ref{cropro-auto}. But since these constructions will not be needed, we do not go into the details here.

\subsection{Examples of subsemigroups}

It is not clear for which semigroups $\pi_s: C^*(P) \to C^*_s(P)$ is an isomorphism. But in typical examples, we see that condition III$_G$ is already satisfied in $C^*(P)$.

For instance, let $(G,P)$ be a quasi-lattice ordered group as in \S~\ref{exp}. In that case, III$_G$ is automatically satisfied in $C^*(P)$. Namely, given $p$, $q$ in $P$ such that $(pP) \cap (qP) \neq \emptyset$, we can find $r \in P$ such that $(pP) \cap (qP) = rP$, and then $v_p^* v_q = v_p^* v_p v_p^* v_q v_q^* v_q = v_p^* v_r v_r^* v_q = v_{p^{-1} r} v_{q^{-1} r}^*$. Applying this several times, we can write $v_{p_1}^* v_{q_1} \dotsm v_{p_m}^* v_{q_m}$ as $v_x v_y^*$ for some $x, y \in P$ if $\e{q_m^{-1} p_m \dotsm q_1^{-1} p_1 P} \neq 0$. Now if $p_1^{-1} q_1 \dotsm p_m^{-1} q_m = e$ in $G$, then $xy^{-1} = e$ in $G$, hence $x=y$. Therefore, $v_{p_1}^* v_{q_1} \dotsm v_{p_m}^* v_{q_m} = v_x v_x^*$ is a projection, and we deduce $v_{p_1}^* v_{q_1} \dotsm v_{p_m}^* v_{q_m} = (v_{p_1}^* v_{q_1} \dotsm v_{p_m}^* v_{q_m})^* (v_{p_1}^* v_{q_1} \dotsm v_{p_m}^* v_{q_m}) = \e{q_m^{-1} p_m \dotsm q_1^{-1} p_1 P}$ in $C^*(P)$.

Another class of such examples is given by left Ore semigroups.

\bdefin
A semigroup $P$ is called right reversible if for every $p$, $q$ in $P$, we have $(Pp) \cap (Pq) \neq \emptyset$.
\edefin

\bdefin
A semigroup is called left Ore if it is cancellative (i.e. left and right cancellative) and right reversible.
\edefin

We have the following
\btheo[Ore, Dubreil]
\label{OD}
A semigroup $P$ can be embedded into a group $G$ such that $G = P^{-1} P = \menge{q^{-1} p}{p,q \in P}$ if and only if $P$ is left Ore.
\etheo
The reader may consult \cite{Cl-Pr}, Theorem~1.24 or \cite{La}, \S~1.1 for more explanations about this theorem. For later purposes, we also introduce the following
\bdefin
A semigroup $P$ is called left reversible if for every $p$, $q$ in $P$, we have $(pP) \cap (qP) \neq \emptyset$.
\edefin

\bdefin
A semigroup is called right Ore if it is cancellative and left reversible.
\edefin

The analogue of Theorem~\ref{OD} is
\btheo[Ore, Dubreil (right version)]
\label{OD-rightquot}
A semigroup $P$ can be embedded into a group $G$ such that $G = P P^{-1} = \menge{p q^{-1}}{p,q \in P}$ if and only if $P$ is right Ore.
\etheo

Now let us see that for a left Ore semigroup, condition III$_G$ is already satisfied in $C^*(P)$. Given $p$, $q$ in $P$, there exist by right reversibility $r$, $s$ in $P$ such that $rp = sq$. Thus $v_q v_p^* = v_s^* v_s v_q v_p^* = v_s^* v_r v_p v_p^*$. Applying this several times, we can write $v_{p_1}^* v_{q_1} \dotsm v_{p_m}^* v_{q_m}$ as $v_y^* v_x e_X$ for some $X \in \cJ$. If $p_1^{-1} q_1 \dotsm p_m^{-1} q_m = e$ holds in $G = P^{-1} P$, then $y^{-1}x = e$ in $G$, hence $x=y$. Thus we again conclude that $v_{p_1}^* v_{q_1} \dotsm v_{p_m}^* v_{q_m} = v_x^* v_x e_X$ is a projection, and the same argument as in the quasi-lattice ordered case gives $v_{p_1}^* v_{q_1} \dotsm v_{p_m}^* v_{q_m} = \e{q_m^{-1} p_m \dotsm q_1^{-1} p_1 P}$ in $C^*(P)$.

\subsection{Conditional expectations}

We conclude this section with a few observations which will be used later on. First of all, there is a faithful conditional expectation $\cE_r: \cL(\ell^2(P)) \to \ell^{\infty}(P) \subseteq \cL(\ell^2(P))$ characterized by
\bgloz
  \spkl{\cE_r(T) \varepsilon_x, \varepsilon_x} = \spkl{T \varepsilon_x, \varepsilon_x} \fa T \in \cL(\ell^2(P)), x \in P.
\egloz
Here $\ell^\infty(P)$ acts on $\ell^2(P)$ by multiplication operators.
\blemma
\label{fce}
If $P$ embeds into a group $G$, then $\cE_r(C^*_r(P)) = D_r(P)$.
\elemma
\bproof
As $D_r(P) \subseteq \ell^{\infty}(P)$, it is clear that $\cE_r(C^*_r(P))$ contains $D_r(P)$. It remains to prove \an{$\subseteq$}. By the definition of the reduced semigroup C*-algebra, we have
\bgloz
  C^*_r(P) = \clspan(\menge{V_{p_1}^* V_{q_1} \dotsm V_{p_m}^* V_{q_m}}{m \in \Zz \pos \text{; } p_i, q_i \in P \fa 1 \leq i \leq m}).
\egloz
So it suffices to prove that for every $p_1, q_1, \dotsc, p_m, q_m \in P$, $\cE_r(V_{p_1}^* V_{q_1} \dotsm V_{p_m}^* V_{q_m}) \in D_r(P)$. Set $V \defeq V_{p_1}^* V_{q_1} \dotsm V_{p_m}^* V_{q_m}$. It is clear that for every $x \in P$, $V \varepsilon_x$ is either $0$ or of the form $\varepsilon_y$ for some $y \in P$. Now assume that $\cE_r(V) \neq 0$. Then there must be $x \in P$ with $V \varepsilon_x = \varepsilon_x$. But this implies that $p_1^{-1} q_1 \dotsm p_m^{-1} q_m x = x$, and thus $p_1^{-1} q_1 \dotsm p_m^{-1} q_m = e$ in $G$. Lemma~\ref{iii_G-red} implies that $V = \E{q_m^{-1} p_m \dotsm q_1^{-1} p_1 P}$ lies in $D_r(P)$.
\eproof
\bremark
This lemma implies that $D_r(P) = C^*_r(P) \cap \ell^\infty(P)$ if $P$ embeds into a group. At this point, we see that it is convenient to work the the family $\cJ$ which is closed under pre-images (with respect to left multiplication), see Remark~\ref{cJ'}.
\eremark

Now let $P$ be a subsemigroup of a group $G$, and let $\stalg(P)$ be the sub-*-algebra of $C^*_s(P)$ generated by the $v_p$, $p \in P$. Set for $g \in G$
\bgl
\label{def-D_g}
  D_g \defeq \lspan(\menge{v_{p_1}^* v_{q_1} \dotsm v_{p_m}^* v_{q_m}}{m \geq 1 \text{; } p_i, q_i \in P \text{ and } p_1^{-1} q_1 \dotsc p_m^{-1} q_m = g})
\egl
as a subspace of $\stalg(P)$. We then obviously have $\stalg(P) = \sum_{g \in G} D_g$.

\blemma
\label{E}
Assume that $P$ embeds into a group $G$ and that the constructible right ideals of $P$ are independent. Then there is a conditional expectation $\cE_s: C^*_s(P) \to D_s(P)$ with
\bgln
\label{E/D}
  && \cE_s \vert_{D_g} = 0 \falls g \neq e \text{ and } \cE_s \vert_{D_e} = \id_{D_e}; \\
\label{kerlambda=kerE}
  && \ker(\lambda) \cap {C^*_s(P)}_+ = \ker(\cE_s) \cap {C^*_s(P)}_+,
\egln
where ${C^*_s(P)}_+$ denotes the set of positive elements in $C^*_s(P)$.
\elemma
\bproof
Since we assume that the constructible right ideals of $P$ are independent, we know that $\lambda \vert_{D_s(P)}$ is an isomorphism. Thus we can set
\bgloz
  \cE_s \defeq (\lambda \vert_{D_s(P)})^{-1} \circ \cE_r \circ \lambda: C^*_s(P) \to D_s(P).
\egloz
We have
\bgloz
  \cE_r(V_{p_1}^* V_{q_1} \dotsm V_{p_m}^* V_{q_m}) =
  \bfa
  \E{q_m^{-1} p_m \dotsm q_1^{-1} p_1 P} \text{ if } p_1^{-1} q_1 \dotsc p_m^{-1} q_m = e, \\
  0 \text{ if } p_1^{-1} q_1 \dotsc p_m^{-1} q_m \neq e. 
  \efa
\egloz
Therefore we obviously have $\cE_s \vert_{D_g} = 0$ if $g \neq e$. And for $p_1, q_1, \dotsc, p_m, q_m \in P$ with $p_1^{-1} q_1 \dotsm p_m^{-1} q_m = e$ in $G$, we have
\bglnoz
  && \cE_s (v_{p_1}^* v_{q_1} \dotsm v_{p_m}^* v_{q_m}) = ((\lambda \vert_{D_s(P)})^{-1} \circ \cE_r)(V_{p_1}^* V_{q_1} \dotsm V_{p_m}^* V_{q_m}) \\
  &=& (\lambda \vert_{D_s(P)})^{-1} (\E{q_m^{-1} p_m \dotsm q_1^{-1} p_1 P}) = \e{q_m^{-1} p_m \dotsm q_1^{-1} p_1 P} 
  \overset{\text{III}_G}{=} v_{p_1}^* v_{q_1} \dotsm v_{p_m}^* v_{q_m}.
\eglnoz
\eproof

\section{Amenability}
\label{amenability}
In this section, our goal is to study the relationship between semigroups and their semigroup C*-algebras in the context of amenability. It turns out that, using our constructions of semigroup C*-algebras, there are strong parallels between the semigroup case and the group case. Indeed, one of our main goals in this section is to show that the analogues of \cite{Br-Oz}, Chapter~2, Theorem~6.8~(1)--(7) are also equivalent in the case of semigroups (under certain assumptions on the semigroups). Apart from this result, we also prove a few additional statements.

Let us first state our main result. To do so, we recall some definitions. The reader may find more explanations in \cite{Pa}.

\bdefin
A discrete semigroup $P$ is left amenable if there exists a left invariant mean on $\ell^{\infty}(P)$, i.e. a state $\mu$ on $\ell^{\infty}(P)$ such that for every $p \in P$ and $f \in \ell^{\infty}(P)$, $\mu(f(p \sqcup)) = \mu(f)$. Here $f(p \sqcup)$ is the composition of $f$ after left multiplication with $p$.
\edefin
\bdefin
An approximate left invariant mean on a discrete semigroup $P$ is a net $(\mu_i)_i$ in $\ell^1(P)$ of positive elements of norm $1$ with the property that
\bgloz
  \lim_i \norm{\mu_i - \mu_i(p \sqcup)}_{\ell^1(P)} = 0 \fa p \in P.
\egloz
Here $\mu_i(p \sqcup)$ again is the composition of $\mu_i$ after left multiplication with $p$.
\edefin
\bdefin
A discrete semigroup $P$ satisfies the strong F{\o}lner condition if for every finite subset $C \subseteq P$ and every $\varepsilon > 0$, there exists a non-empty finite subset $F \subseteq P$ such that $\abs{(pF) \Delta F} / \abs{F} < \varepsilon \fa p \in C$.
\edefin
Here $\Delta$ stands for symmetric difference.

\subsection{Statements}
\label{stats}
Let $P$ be a discrete left cancellative semigroup. We consider the following statements:
\begin{itemize}
\item[1)] $P$ is left amenable.
\item[2)] $P$ has an approximate left invariant mean.
\item[3)] $P$ satisfies the strong F{\o}lner condition.
\item[4)] There exists a net $(\xi_i)_i$ in $\ell^2(P)$ with $\norm{\xi_i} = 1$ for all $i$ and $\lim_i \norm{V_p \xi_i - \xi_i} = 0$ for all $p \in P$.
\item[5)] There exists a net $(\xi_i)_i$ in $C_c(P) \subseteq \ell^2(P)$ with $\norm{\xi_i} = 1$ for all $i$ such that $\lim_i \spkl{V_{p_1}^* V_{q_1} \dotsm V_{p_n}^* V_{q_n} \xi_i, \xi_i} = 1$ for all $n \in \Zz \pos$; $p_1, q_1, \dotsc, p_n, q_n \in P$.
\item[6)] The left regular representation $\lambda: C^*_s(P) \to C^*_r(P)$ is an isomorphism and there exists a non-zero character on $C^*_s(P)$.
\item[7)] There exists a non-zero character on $C^*_r(P)$.
\end{itemize}
Our goal is to show that for a discrete left cancellative semigroup, we always have \an{1) $\LRarr$ 2) $\LRarr$ 3) $\Rarr$ 4) $\Rarr$ 5)} and \an{6) $\Rarr$ 7) $\Rarr$ 1)}, and that if $P$ is also right cancellative and if the constructible right ideals are independent (see Definition~\ref{ind}), then \an{5) $\Rarr$ 6)} holds as well. With Corollary~\ref{un-cond->lrr/D} in mind, it is not surprising that independence of the family of constructible right ideals plays a role in the context of amenability. Moreover, note that 6) only makes sense if $P$ can be embedded into a group. Thus our assumption that $P$ should be cancellative is certainly necessary, and as a part of \an{5) $\Rarr$ 6)}, we will prove that 5) implies that $P$ embeds into a group. In addition, we will see in Remark~\ref{5->7} that 5) implies 7) for every discrete left cancellative semigroup.

Before we start with the proofs, let us remark that the equivalence of 1), 2) and 3) for discrete left cancellative semigroups is certainly known, and that these equivalences can be proven as in the group case. We include proofs of these equivalences for the sake of completeness. Moreover, the implications \an{3) $\Rarr$ 4) $\Rarr$ 5)} and \an{6) $\Rarr$ 7)} are easy. And for the implication \an{7) $\Rarr$ 1)}, the proof in the group case as presented in \cite{Br-Oz}, Chapter~2, Theorem~6.8 carries over to the case of semigroups. Again, for the sake of completeness, we present a proof for this implication. Both for the equivalence of 1), 2) and 3) as well as for the implication \an{7) $\Rarr$ 1)}, we only have to check that in the proofs of the corresponding statements in the group case, we can avoid taking inverses as this is in general not possible in semigroups. And finally, to prove \an{5) $\Rarr$ 6)} under the additional assumptions that $P$ is right cancellative and that the constructible right ideals of $P$ are independent, we adapt A. Nica's ideas in \cite{Ni}, \S~4.4 to our situation.

\subsection{Proofs}
We start with {\bf \an{1) $\LRarr$ 2)}}. First assume that there is a left invariant mean $\mu$ on $\ell^{\infty}(P)$. As the unit ball of $\ell^1(P)$ is weak*-dense in the unit ball of $\ell^1(P)'' \cong \ell^{\infty}(P)'$, there exists a net $(\mu_i)_i$ of positive elements in $\ell^1(P)$ with norm $1$ which converges to $\mu$ in the weak*-topology. This means that $\lim_i \mu_i(f) = \mu(f)$ for every $f \in \ell^{\infty}(P)$. We want to show that for every $p \in P$ and $f \in \ell^{\infty}(P)$, $\lim_i \mu_i(f) - (\mu_i(p \sqcup))(f) = 0$. To prove this, take $f \in \ell^{\infty}(P)$, $p \in P$ and define a function $g \in \ell^{\infty}(P)$ by
$
g(q) \defeq 
  \bfa 
    f(r) \falls q=pr \\
    0 \sonst.
  \efa
$
Then $\lim_i (\mu_i(g(p \sqcup)) - \mu_i(g)) = \mu(g(p \sqcup)) - \mu(g) = 0$ as $\mu$ is left invariant. At the same time,
\bglnoz
  && \mu_i(g(p \sqcup)) - \mu_i(g) = \sum_q \mu_i(q) g(pq) - \sum_q \mu_i(q) g(q) \\
  &=& \sum_q \mu_i(q) g(pq) - \sum_q \mu_i(pq) g(pq) - \sum_{q \notin pP} \mu_i(q) \underbrace{g(q)}_{=0} \\
  &=& \sum_q \mu_i(q) f(q) - \sum_q \mu_i(pq) f(q) = \mu_i(f) - (\mu_i(p \sqcup))(f).
\eglnoz
This shows that we indeed have $\lim_i \mu_i(f) - (\mu_i(p \sqcup))(f) = 0$. Hence, for every $n \in \Zz \pos$ and $p_1, \dotsc, p_n \in P$, $(0, \dotsc, 0)$ lies in the weak closure of
\bgl
\label{1-2}
  \menge{(\nu - \nu(p_j \sqcup))_{j=1, \dotsc, n}}{\nu \in \ell^1(P), \nu \geq 0, \norm{\nu} \leq 1}.
\egl
As this set is convex, it follows from the Hahn-Banach separation theorem that its weak and norm closures coincide. That $(0, \dotsc, 0)$ lies in the norm closure of \eqref{1-2} tells us that $P$ has an approximate left invariant mean. This proves \an{1) $\Rarr$ 2)}.

For the reverse implication, assume that $P$ has an approximate left invariant mean $(\mu_i)_i$. By definition, this means
\bgl
\label{2-1_1}
  \lim_i \norm{\mu_i - \mu_i(p \sqcup)}_{\ell^1(P)} = 0 \fa p \in P.
\egl
Moreover, we have $\norm{\mu_i - \mu_i(p \sqcup)}_{\ell^1(P)} \geq \norm{\mu_i}_{\ell^1(P)} - \norm{\mu_i(p \sqcup)}_{\ell^1(P)} = \sum_{q \notin pP} \abs{\mu_i(q)}$. It follows that
\bgl
\label{2-1_2}
  \lim_i \sum_{q \notin pP} \abs{\mu_i(q)} = 0.
\egl
Now $\ell^{\infty}(P)' \cong \ell^1(P)''$, and by the theorem of Banach-Alaoglu, the unit ball of $\ell^1(P)''$ is weak*-compact. Hence by passing to a suitable subnet if necessary, we may assume that the net $(\mu_i)_i$ converges to an element $\mu \in \ell^1(P)'' \cong \ell^{\infty}(P)'$ in the weak*-topology. $\mu$ has to be a state on $\ell^{\infty}(P)$ as the $\mu_i$ are positive with norm $1$. For every $f \in \ell^{\infty}(P)$ and $p \in P$ we have
\bglnoz
  && \abs{\mu(f(p \sqcup)) - \mu(f)} = \lim_i \abs{\mu_i(f(p \sqcup)) - \mu_i(f)} \\
  &=& \lim_i \left\vert \sum_{q \in P} \mu_i(q) f(pq) - \sum_{q \in P} \mu_i(q) f(q) \right\vert \\
  &=& \lim_i \left\vert \sum_{q \in P} (\mu_i(q)-\mu_i(pq)) f(pq) - \sum_{q \notin pP} \mu_i(q) f(q) \right\vert \\
  &\leq& \lim_i 
  \rukl{\norm{\mu_i - \mu_i(p \sqcup)}_{\ell^1(P)} \cdot \norm{f}_{\ell^{\infty}(P)} + \sum_{q \notin pP} \abs{\mu_i(q)} \norm{f}_{\ell^{\infty}(P)}} \\
  &=& 0 
\eglnoz
by \eqref{2-1_1} and \eqref{2-1_2}. Thus $\mu$ is a left invariant mean. This proves \an{2) $\Rarr$ 1)}.

Let us prove {\bf \an{1) $\LRarr$ 3)}}. First of all, if $P$ has an approximate left invariant mean $(\mu_i)_i$, then we always have
\bgl
\label{lim-with-inv}
  \lim_i \norm{\mu_i(p^{-1} \sqcup) - \mu_i}_{\ell^1(P)} = 0,
\egl
where 
$
  \mu_i(p^{-1} \sqcup) (q) = 
  \bfa
  \mu_i(q') \falls q = pq' \text{ for some } q' \in P \\
  0 \falls q \notin pP
  \efa
$. The reason is that we have
\bglnoz
  && \norm{\mu_i(p^{-1} \sqcup) - \mu_i}_{\ell^1(P)} = \sum_{q \in pP} \abs{\mu_i(p^{-1} \sqcup)(q) - \mu_i(q)} + \sum_{q \notin pP} \abs{\mu_i(q)} \\
  &=& \sum_{q' \in P} \abs{\mu_i(q') - \mu_i(pq')} + \sum_{q \notin pP} \abs{\mu_i(q)} = \norm{\mu_i - \mu_i(p \sqcup)}_{\ell^1(P)} + \sum_{q \notin pP} \abs{\mu_i(q)}
\eglnoz
and $\lim_i \sum_{q \notin pP} \abs{\mu_i(q)} = 0$ by \eqref{2-1_2}.

Now, assume that $P$ has an approximate left invariant mean. Let $C$ be a finite subset $P$ and let $\varepsilon > 0$ be given. By 2) and the fact proven above that every approximate left invariant mean $(\mu_i)_i$ satisfies \eqref{lim-with-inv}, there exists a positive $\ell^1$-function $\mu$ of $\ell^1$-norm $1$ with
\bgl
\label{2-3}
  \sum_{p \in C} \norm{\mu(p^{-1} \sqcup) - \mu}_{\ell^1(P)} < \varepsilon.
\egl
For $t \in [0,1]$, we set $F(\mu,t) \defeq \menge{q \in P}{\mu(q)>t}$. We claim that for a suitable choice of $t$, the inequality $\max_{p \in C} \abs{pF(\mu,t) \Delta F(\mu,t)} / \abs{F(\mu,t)} < \varepsilon$ holds. We have
\bglnoz
  && \norm{\mu(p^{-1} \sqcup) - \mu}_{\ell^1(P)} = \sum_{q \in P} \abs{(\mu(p^{-1} \sqcup) - \mu)(q)} \\
  &=& \sum_{q \in P} \int_0^1 \abs{\1z_{[0,\mu(p^{-1} \sqcup)(q)]}(t) - \1z_{[0,\mu(q)]}(t)} dt \\
  &=& \sum_{q \in P} \int_0^1 \abs{\1z_{F(\mu(p^{-1} \sqcup),t)}(q) - \1z_{F(\mu,t)}(q)} dt = \int_0^1 \abs{(pF(\mu,t)) \Delta F(\mu,t)} dt
\eglnoz
and
\bglnoz
  && \int_0^1 \varepsilon \abs{F(\mu,t)} dt = \varepsilon \int_0^1 \sum_{q \in P} \1z_{F(\mu,t)}(q) dt = \varepsilon \sum_{q \in P} \int_0^1 \1z_{F(\mu,t)}(q) dt \\
  &=& \varepsilon \sum_{q \in P} \int_0^1 \1z_{[0,\mu(q)]}(t) dt = \varepsilon \sum_{q \in P} \mu(q) = \varepsilon.
\eglnoz
Plugging these two inequalities into \eqref{2-3}, we obtain
\bgloz
  \int_0^1 \varepsilon \abs{F(\mu,t)} dt > \int_0^1 \sum_{p \in C} \abs{(pF(\mu,t)) \Delta F(\mu,t)} dt
\egloz
Thus there is $t \in [0,1]$ with $\varepsilon \abs{F(\mu,t)} > \sum_{p \in C} \abs{(pF(\mu,t)) \Delta F(\mu,t)}$. Therefore $P$ satisfies the strong F{\o}lner condition. So we have proven \an{2) $\Rarr$ 3)}.

To prove the reverse implication, observe that 3) tells us that there exists a net $(F_i)_i$ of non-empty finite subsets of $P$ such that $\lim_i \abs{(pF_i) \Delta F_i} / \abs{F_i} = 0 \fa p \in P$. Set $\mu_i \defeq \tfrac{1}{\abs{F_i}} \1z_{F_i}$. It is clear that $(\mu_i)_i$ is a net of positive $\ell^1$-functions of $\ell^1$-norm $1$. Moreover, $\norm{\mu_i - \mu_i(p \sqcup)}_{\ell^1(P)} \leq \norm{\mu_i(p^{-1} \sqcup) - \mu_i}_{\ell^1(P)} = \abs{\tfrac{1}{\abs{F_i}}(\1z_{pF_i} - \1z_{F_i})}_{\ell^1(P)} = \abs{(pF_i) \Delta F_i} / \abs{F_i} \lori_i 0$ for all $p$ in $P$. Thus $(\mu_i)_i$ is an approximate left invariant mean. This proves \an{3) $\Rarr$ 2)}.

To prove {\bf \an{3) $\Rarr$ 4)}}, first note that since $P$ satisfies the strong F{\o}lner condition, there is a net $(F_i)_i$ of non-empty finite subsets of $P$ with $\lim_i \abs{(pF_i) \Delta F_i} / \abs{F_i} = 0$ for all $p \in P$. Now set $\xi_i \defeq \abs{F_i}^{-\halb} \1z_{F_i}$. Here $\1z_{F_i}$ is the characteristic function of $F_i \subseteq P$. It is clear that every $\xi_i$ lies in $\ell^2(P)$ and has norm $1$. Moreover, for every $p \in P$, $V_p \xi_i - \xi_i = \abs{F_i}^{-\halb} (\1z_{pF_i} - \1z_{F_i})$. It follows that $\norm{V_p \xi_i - \xi_i}^2 = \abs{(pF_i) \Delta F_i} / \abs{F_i} \lori_i 0 \fa p \in P$. This proves \an{3) $\Rarr$ 4)}.

{\bf \an{4) $\Rarr$ 5)}:} By an approximation argument, we can without loss of generality assume that the $\xi_i$ from 4) all lie in $C_c(P)$. We have by 4) that $\lim_i \norm{V_p \xi_i - \xi_i} = 0$ for all $p \in P$ and also $\norm{V_p^* \xi_i - \xi_i} \leq \norm{V_p^*} \cdot \norm{\xi_i - V_p \xi_i} \lori_i 0$ for all $p \in P$. Hence
\bglnoz
  && \abs{\spkl{V_{p_1}^* V_{q_1} \dotsm V_{p_n}^* V_{q_n} \xi_i, \xi_i} - 1} \\
  &=& \left\vert
  \sum_{j=1}^n 
  \left(
  \spkl{V_{p_1}^* V_{q_1} \dotsm V_{p_j}^* V_{q_j} \xi_i, \xi_i} 
  - \spkl{V_{p_1}^* V_{q_1} \dotsm V_{p_{j-1}}^* V_{q_{j-1}} V_{p_j}^* \xi_i, \xi_i} \right. \right. \\
  &&
  \left. \left.
  + \spkl{V_{p_1}^* V_{q_1} \dotsm V_{p_{j-1}}^* V_{q_{j-1}} V_{p_j}^* \xi_i, \xi_i}
  - \spkl{V_{p_1}^* V_{q_1} \dotsm V_{p_{j-1}}^* V_{q_{j-1}} \xi_i, \xi_i}
  \right)
  \right\vert \\
  &\leq& \sum_{j=1}^n \norm{V_{q_j} \xi_i - \xi_i} + \norm{V_{p_j}^* \xi_i - \xi_i} \lori_i 0
\eglnoz
for all $n \in \Zz \pos$ and $p_1, q_1, \dotsc, p_n, q_n \in P$. This proves \an{4) $\Rarr$ 5)}.

{\bf \an{6) $\Rarr$ 7)}} is trivial.

For {\bf \an{7) $\Rarr$ 1)}}, let $\chi: C^*_r(P) \to \Cz$ be a non-zero character. Viewing $\chi$ as a state, we can extend it by the theorem of Hahn-Banach to a state on $\cL(\ell^2(P))$. We then restrict the extension to $\ell^{\infty}(P) \subseteq \cL(\ell^2(P))$ and call this restriction $\mu$. The point is that by construction, $\mu \vert_{C^*_r(P)} = \chi$ is multiplicative, hence $C^*_r(P)$ is in the multiplicative domain of $\mu$. Thus we obtain for every $f \in \ell^{\infty}(P)$ and $p \in P$
\bgloz
  \mu(f(p \sqcup)) = \mu(V_p^* f V_p) = \mu(V_p^*) \mu(f) \mu(V_p) = \mu(V_p)^* \mu(V_p) \mu(f) = \mu(f).
\egloz
Thus $\mu$ is a left invariant mean on $\ell^{\infty}(P)$. Hence we have proven \an{7) $\Rarr$ 1)}.

It remains to discuss the implication {\bf \an{5) $\Rarr$ 6)}}. We start with the following
\blemma
\label{5-left-rev}
5) implies that $P$ is left reversible.
\elemma
\bproof
Let $(\xi_i)_i$ be a net as in 5). For $p_1, p_2 \in P$, we have $\lim_i \spkl{V_{p_1} V_{p_1}^* V_{p_2} V_{p_2}^* \xi_i, \xi_i} = 1$. In particular, $V_{p_1} V_{p_1}^* V_{p_2} V_{p_2}^* \neq 0$. But $V_{p_1} V_{p_1}^* V_{p_2} V_{p_2}^* = \E{(p_1 P) \cap (p_2 P)}$, hence $(p_1 P) \cap (p_2 P) \neq \emptyset$. This shows that $P$ is left reversible.
\eproof
\bcor
\label{5-embed}
If $P$ is cancellative and 5) holds, then $P$ embeds into a group $G$ such that $G = P P^{-1}$.
\ecor
\bproof
This follows from the previous lemma and Theorem~\ref{OD-rightquot}.
\eproof

\blemma
\label{left-rev-char}
A subsemigroup $P$ of a group is left reversible if and only if there exists a non-zero character on $C^*_s(P)$.
\elemma
\bproof
If $\chi$ is a non-zero character on $C^*_s(P)$, then for every $p_1, p_2 \in P$, we have $\chi(\e{(p_1 P) \cap (p_2 P)}) = \chi(v_{p_1} v_{p_1}^* v_{p_2} v_{p_2}^*) = \chi(v_{p_1}) \chi(v_{p_1}^*) \chi(v_{p_2}) \chi(v_{p_2}^*) = 1$. This implies that $(p_1 P) \cap (p_2 P) \neq \emptyset$ because otherwise $\e{(p_1 P) \cap (p_2 P)}$ would vanish.

If $P$ is left reversible, then by universal property of $C^*_s(P)$, there is a homomorphism $C^*_s(P) \to \Cz$ sending $C^*_s(P) \ni v_p$ to $1 \in \Cz$ and $C^*_s(P) \ni e_X$ to $1 \in \Cz$ if $X \neq \emptyset$ and to $0 \in \Cz$ if $X = \emptyset$ for every $p \in P$ and $X \in \cJ$. This is compatible with relation III$_G$ as $q_m^{-1} p_m \dotsm q_1^{-1} p_1 P$ is never empty. The last fact follows inductively on $m$ using the observation that for every non-empty right ideal $X$ of $P$, we have $q^{-1} p X = q^{-1} ((pX) \cap (qP))$, and that $(pX) \cap (qP) \neq \emptyset$ by left reversibility.
\eproof

It remains to prove that 5) implies that $\lambda: C^*_s(P) \to C^*_r(P)$ is an isomorphism if $P$ is cancellative (not only left cancellative, but also right cancellative) and if the constructible right ideals of $P$ are independent. Recall the definition of $D_g$ from \eqref{def-D_g}. For a positive functional $\varphi$ on $C^*_s(P)$, we define the d-support of $\varphi$ as $\dsupp(\varphi) \defeq \menge{g \in G}{\varphi \vert_{D_g} \neq 0}$. Moreover, we set
\bgloz
  \cV \defeq \menge{v_{p_1}^* v_{q_1} \dotsm v_{p_n}^* v_{q_n}}{n \in \Zz \pos; p_i, q_i \in P}.
\egloz

Our aim is to show
\btheo
\label{thm}
Let $P$ be a subsemigroup of a group $G$, and assume that the constructible right ideals of $P$ are independent. If there exists a net $(\varphi_i)_i$ of states on $C^*_s(P)$ with finite d-support such that $\lim_i \varphi_i(v) = 1$ for every $0 \neq v$ in $\cV$, then $\lambda: C^*_s(P) \to C^*_r(P)$ is an isomorphism.
\etheo
Note that this is the analogue of the implication \an{(5) $\Rarr$ (6)} in \cite{Br-Oz}, Chapter~2, Theorem~6.8 in the group case. To prove the theorem, we first show
\blemma
\label{lem-ineq}
Let $\varphi$ be a positive functional on $C^*_s(P)$ with finite d-support. We then have for all $x \in C^*_s(P)$:
\bgl
\label{5-6}
  \abs{\varphi(x)}^2 \leq \abs{\dsupp(\varphi)} \norm{\varphi} \varphi(\cE_s(x^*x)).
\egl
\elemma
Here $\cE_s$ is the conditional expectation from Lemma~\ref{E}.
\bproof
It certainly suffices to prove our assertion for $x$ in $\stalg(P) = \sum_{g \in G} D_g$. Take such an element $x$. Let $\dsupp(\varphi) = \gekl{g_1, \dotsc, g_n}$. We can find a finite subset $F \subseteq G$ so that $x = \sum_{g \in F} x_g \text{ with } x_g \in D_g$ and $\dsupp(\varphi) \subseteq F$, i.e. $\gekl{g_1, \dotsc, g_n} \subseteq F$. Then $\varphi(x) = \sum_{g \in F} \varphi(x_g) = \sum_{j=1}^n \varphi(x_{g_j})$. Thus, using the Cauchy-Schwarz inequality twice, we obtain
\bglnoz
  && \abs{\varphi(x)}^2 = \left\vert \sum_{j=1}^n \varphi(x_{g_j}) \right\vert^2 = \abs{\spkl{(\varphi(x_{g_j}))_j,(1)_j}_{\Cz^n}}^2
  \leq \norm{(\varphi(x_{g_j}))_j}_{\Cz^n}^2 \norm{(1)_j}_{\Cz^n}^2 \\
  &=& n \sum_{j=1}^n \abs{\varphi(x_{g_j})}^2 
  = n \sum_{j=1}^n \abs{\spkl{x_{g_j},1}_{\varphi}}^2 \leq n \norm{\varphi} \sum_{j=1}^n \varphi(x_{g_j}^* x_{g_j}).
\eglnoz
Hence it suffices to prove $\sum_{j=1}^n x_{g_j}^* x_{g_j} \leq \cE_s(x^* x)$. We have by \eqref{E/D} and because of $D_g^* D_h \subseteq D_{g^{-1}h}$ for all $g,h \in G$ that
\bgloz
  \cE_s(x^* x) = \sum_{g,h \in F} \cE_s(x_g^* x_h) = \sum_{g,h \in F} \delta_{g,h} x_g^* x_h = \sum_{g \in F} x_g^* x_g \geq \sum_{j=1}^n x_{g_j}^* x_{g_j}.
\egloz
This proves \eqref{lem-ineq}.
\eproof

\bprop
\label{finite-supp-dense}
$\lambda: C^*_s(P) \to C^*_r(P)$ is an isomorphism if the set of positive functionals on $C^*_s(P)$ with finite d-support is dense in the space of all positive functionals on $C^*_s(P)$ in the weak*-topology.
\eprop
\bproof
Take $x \in \ker(\lambda)$. Passing over to $x^* x$ if necessary, we may assume $x \geq 0$. Take a positive functional $\varphi$ on $C^*_s(P)$ with finite d-support. We then have because of $\lambda(x)=0$ that $\lambda(x^* x)=0$, thus $\cE_s(x^* x)=0$ by \eqref{kerlambda=kerE}. Hence it follows from \eqref{lem-ineq} that $\varphi(x)=0$. So we have shown that $\varphi(x)=0$ for every positive functional on $C^*_s(P)$ with finite d-support. By our assumption in the proposition, the positive functionals with finite d-support are weak*-dense in the space of all positive functionals. Hence $\varphi(x)=0$ for every positive functional $\varphi$ on $C^*_s(P)$. This however implies that $x=0$. We conclude that $\lambda$ must be injective, hence an isomorphism.
\eproof
Actually, the converse of the proposition is valid as well, and is simpler to prove. To proceed, we need another
\blemma
\label{pos-fct}
Let $\varphi$ and $\phi$ be positive functionals on $C^*_s(P)$. Then there exists a unique positive functional $\psi$ on $C^*_s(P)$ such that $\psi(v) = \varphi(v) \phi(v)$ for all $v \in \cV$.
\elemma
\bproof
Just set $\psi = (\varphi \otimes \phi) \circ \Delta$ with $\Delta$ given by \eqref{Delta-C_s}.
\eproof
Finally, with all these preparations, we can prove our theorem.
\bproof[Proof of Theorem~\ref{thm}]
Let $\phi$ be a positive functional on $C^*_s(P)$. Let $\varphi_i$ be the states given by the hypothesis of our theorem, they satisfy
\bgl
\label{lim-phi}
  \lim_i \varphi_i(v) = 1 \text{ for every } 0 \neq v \in \cV.
\egl
By Lemma~\ref{pos-fct}, there exists a net $(\phi_i)_i$ of positive functionals on $C^*_s(P)$ such that for all $i$, 
\bgl
\label{phi-phiphi}
  \phi_i(v) = \varphi_i(v) \phi(v) \fa v \in \cV.
\egl
In particular, $\norm{\phi_i} = \norm{\phi}$ since $\phi_i(1) = \phi(1) = \norm{\phi}$. It is then clear that for every $i$, $\dsupp(\phi_i) \subseteq \dsupp(\varphi_i)$ is finite. Moreover, we have $\lim_i \phi_i(v) = \phi(v) \fa v \in \cV$. This is clear if $v = 0$, and if $v \neq 0$ it follows from \eqref{phi-phiphi} and \eqref{lim-phi}. Thus $\lim_i \phi_i(x) = \phi(x)$ for all $x \in \stalg(P)$, and since $\norm{\phi_i} = \norm{\phi}$ for all $i$, we conclude that we actually have $\lim_i \phi_i(x) = \phi(x)$ for all $x \in C^*_s(P)$. In other words, the net $(\phi_i)_i$ converges to $\phi$ in the weak*-topology. Thus we have seen that the positive functionals with finite d-support are weak*-dense in the space of all positive functionals. By Proposition~\ref{finite-supp-dense}, this implies that $\lambda: C^*_s(P) \to C^*_r(P)$ is an isomorphism. This completes the proof of our theorem.
\eproof

{\bf \an{5) $\Rarr$ 6)} if $P$ is cancellative and if the constructible right ideals of $P$ are independent:} Assume that $P$ is cancellative and that the constructible right ideals of $P$ are independent. We have already seen that 5) implies that $P$ is left reversible in Lemma~\ref{5-left-rev}. Thus $P$ embeds into a group by Corollary~\ref{5-embed}, and there is a non-zero character on $C^*_s(P)$ by Lemma~\ref{left-rev-char}. It remains to prove that $\lambda: C^*_s(P) \to C^*_r(P)$ is an isomorphism. By Theorem~\ref{thm}, it suffices to prove that there exists a net $(\varphi_i)_i$ of states on $C^*_s(P)$ with finite d-support such that $\lim_i \varphi_i(v) = 1$ for every $0 \neq v \in \cV$.

Now take the net $(\xi_i)_i$ in $C_c(P)$ from 5), and set for all $i$: $\varphi_i(x) \defeq \spkl{\lambda(x) \xi_i,\xi_i}$ for every $x \in C^*_s(P)$. It is clear that these $\varphi_i$ are states and that we have $\lim_i \varphi_i(v) = 1$ for every $0 \neq v \in \cV$. Moreover, for every $i$, set $\text{supp}(\xi_i) \defeq \menge{p \in P}{\xi(p) \neq 0}$. By assumption (see 5)), $\text{supp}(\xi_i)$ is a finite set for every $i$. We have $\varphi_i(v_{p_1}^* v_{q_1} \dotsm v_{p_n}^* v_{q_n}) = \spkl{V_{p_1}^* V_{q_1} \dotsm V_{p_n}^* V_{q_n} \xi_i,\xi_i} \neq 0$ only if there are $x$, $y$ in $\text{supp}(\xi_i)$ with $p_1^{-1} q_1 \dotsm p_n^{-1} q_n x = y$. But this implies $p_1^{-1} q_1 \dotsm p_n^{-1} q_n \in (\text{supp}(\xi_i)) (\text{supp}(\xi_i))^{-1}$, or in other words, that $\dsupp(\varphi_i) \subseteq (\text{supp}(\xi_i)) (\text{supp}(\xi_i))^{-1}$. As $\text{supp}(\xi_i)$ is a finite set for every $i$, this proves that for every $i$, $\varphi_i$ has finite d-support. This shows that the conditions in Theorem~\ref{thm} are satisfied, hence that $\lambda: C^*_s(P) \to C^*_r(P)$ is an isomorphism. Thus we have seen that 5) implies 6) if $P$ is cancellative and if the constructible right ideals of $P$ are independent.

\bremark
\label{5->7}
We point out that 5) implies 7) for every discrete left cancellative semigroup $P$. Just set $\chi$ as the weak*-limit of the vector states $\spkl{\sqcup \xi_i,\xi_i}$ of $C^*_r(P)$ where the $\xi_i$ are provided by 5). It is easy to see that $\chi$ is multiplicative.
\eremark

\subsection{Additional results}
There are a few related statements we now turn to. First of all, we can of course consider the following
\bdefin
A discrete semigroup $P$ is called right amenable if there exists a right invariant mean on $\ell^{\infty}(P)$.
\edefin
A right amenable semigroup $P$ is always right reversible, i.e. for every $p_1, p_2 \in P$, we have $(P p_1) \cap (P p_2) \neq \emptyset$. This is the analogue of \cite{Pa}, Proposition~(1.23) if we replace \an{left} in \cite{Pa} by \an{right}. If $P$ is cancellative and right reversible, then $P$ embeds into a group $G$ such that $G = P^{-1} P$ (see Theorem~\ref{OD-rightquot}). $G$ is amenable if $P$ is right amenable (this is the right version of \cite{Pa}, Proposition~(1.27)).

\bprop
\label{P:canc-ra->lambda}
Let $P$ be a cancellative, right amenable semigroup. Then $\lambda\un: {C^*}\un(P) \to C^*_r(P)$ is an isomorphism.
\eprop
\bproof
Consider the embedding $P \into G = P^{-1} P$ from above. We know that ${C^*}\un(P) \cong D\un(P) \rte_{\tau\un} P$ by Lemma~\ref{cropro-descriptions}. By dilation theory for semigroup crossed products by endomorphisms (see \cite{La}), there exists a C*-algebra $D_{\infty}$ with an embedding $D\un(P) \overset{i}{\into} D_{\infty}$ and an action $\tau_{\infty}$ of $G$ on $D_{\infty}$ whose restriction to $P$ leaves $D\un(P)$ invariant and coincides with $\tau\un$. Moreover, $D\un(P) \rte_{\tau\un} P$ embeds into $D_{\infty} \rtimes_{\tau_{\infty}} G$. Let us denote this embedding $D\un(P) \rte_{\tau\un} P \into D_{\infty} \rtimes_{\tau_{\infty}} G$ by $i$ as well.

Since $P$ is right amenable, $G$ is amenable. Hence there is a canonical faithful conditional expectation $\cE_{\infty}$ from $D_{\infty} \rtimes_{\tau_{\infty}} G$ onto $D_{\infty}$. Moreover, using Corollary~\ref{lambda/D-union}, we can construct a conditional expectation on ${C^*}\un(P)$ by setting
\bgl
\label{E-un}
  \cE\un \defeq (\lambda\un \vert_{D\un(P)})^{-1} \circ \cE_r \circ \lambda\un: {C^*}\un(P) \to D\un(P).
\egl
It is easy to see that
$
  \begin{CD}
  D\un(P) \rte_{\tau\un} P @> i >> D_{\infty} \rtimes_{\tau_{\infty}} G \\
  @V \cE\un VV @VV \cE_{\infty} V \\
  D\un(P) @>> i > D_{\infty}
  \end{CD}
$
commutes. But this then shows that $\cE\un$ has to be faithful, and hence that $\lambda\un$ has to be injective (see the Definition of $\cE\un$ in \eqref{E-un}).
\eproof
\bcor
For every cancellative and abelian semigroup $P$, the canonical homomorphism $\lambda\un: {C^*}\un(P) \to C^*_r(P)$ is an isomorphism.
\ecor
\bproof
As remarked in \cite{Pa}, \S~(0.18), every abelian semigroup is amenable.
\eproof
As another consequence of Proposition~\ref{P:canc-ra->lambda}, we obtain an alternative explanation for the result in \cite{C-D-L} that the Toeplitz algebra over the ring of integers $R$ in some number field can be canonically identified with the reduced semigroup C*-algebra of the $ax+b$-semigroup $P_R$ over $R$. First of all, we have proven in Section~\ref{exp} that $\fT[R] \cong C^*(P_R)$. Moreover, we have seen in Lemma~\ref{Dedekind->union} that the constructible right ideals of $P_R$ are independent, so that $\pi\un: C^*(P_R) \to {C^*}\un(P_R)$ is an isomorphism. As $P_R$ embeds into the amenable group $P_K$ (the $ax+b$-group over the quotient field $K$ of $R$) such that $P_K = P_R^{-1} P_R$, it follows that $P_R$ is cancellative, right reversible (see \cite{Cl-Pr}, Theorem~1.24) and hence right amenable (this is the right version of Proposition~(1.28) in \cite{Pa}). Therefore, we may apply Proposition~\ref{P:canc-ra->lambda}. It tells us that $\lambda\un$ is an isomorphism. All in all, we obtain
\bgloz
  \fT[R] \cong C^*(P_R) \overset{\pi\un}{\cong} {C^*}\un(P) \overset{\lambda\un}{\cong} C^*_r(P_R).
\egloz
We point out that the $ax+b$-semigroup over $R$ is not left reversible.

Moreover, we know from the group case that nuclearity of group C*-algebras is closely related to amenability of groups. Here we show
\bprop
Let $P$ be a cancellative, right amenable semigroup. Moreover, assume that $P$ is countable. Then $C^*(P)$, ${C^*}\un(P)$ and $C^*_r(P)$ are nuclear.
\eprop
\bproof
Since we have surjective homomorphisms $C^*(P) \onto {C^*}\un(P) \onto C^*_r(P)$ and because quotients of nuclear C*-algebras are nuclear by \cite{Bla}, Corollary~IV.3.1.13, it suffices to show that our assumptions imply nuclearity of $C^*(P)$.

Using Lemma~\ref{cropro-descriptions} and dilation theory for semigroup crossed products by endomorphisms (see \cite{La}), we conclude that $C^*(P) \cong D(P) \rte_{\tau} P \sim_M D_{\infty} \rtimes_{\tau_{\infty}} G$. Here we use analogous notations as in the proof of Proposition~\ref{P:canc-ra->lambda}. Now $G$ is amenable as $P$ is right amenable, and $D_{\infty}$ is commutative since $D(P)$ is commutative. Hence $D_{\infty} \rtimes_{\tau_{\infty}} G$ is nuclear by \cite{Ror}, Proposition~2.12~(i) and (v). Moreover, all the C*-algebras are separable as $P$ is countable. Hence $C^*(P)$ is nuclear because it is stably isomorphic to a nuclear C*-algebra (see \cite{Ror}, Proposition~2.12~(ii)).
\eproof
In particular, we obtain because every abelian semigroup is amenable:
\bcor
\label{C*(abelian)nuclear}
For every countable, cancellative and abelian semigroup $P$, the C*-algebras $C^*(P)$, ${C^*}\un(P)$ and $C^*_r(P)$ are nuclear.
\ecor

In the reverse direction, we can prove
\bprop
Let $P$ be a cancellative, left reversible semigroup. If $C^*_s(P)$ or ${C^*}\un(P)$ is nuclear, then $P$ is left amenable.
\eprop
\bproof
By assumption, $P$ embeds into a group $G$ with $G = P P^{-1}$ (see Theorem~\ref{OD-rightquot}). As $P$ is left reversible, there exists a canonical projection $C^*_s(P) \to C^*(G)$ sending $v_p$ to $u_p$. Here $u_g$, $g \in G$, denote the unitary generators of $C^*(G)$. As nuclearity passes to quotients by \cite{Bla}, Corollary~IV.3.1.13, nuclearity of $C^*_s(P)$ implies that $C^*(G)$, hence $C^*_r(G)$ must be nuclear as well. By \cite{Br-Oz}, Chapter~2, Theorem~6.8, we conclude that $G$ must be amenable. But a left reversible subsemigroup of an amenable group is itself left amenable by \cite{Pa}, (1.28). The analogous proof works also for ${C^*}\un(P)$ in place of $C^*_s(P)$.
\eproof

\section{Questions and concluding remarks}
An obvious question is: Which semigroups satisfy the condition that their constructible right ideals are independent? It would already be interesting to find out for which integral domains the corresponding $ax+b$-semigroups satisfy this independence condition.

Another question is whether the condition in Lemma~\ref{fce} is actually necessary. In other words, what is the precise relationship between embeddability of $P$ into a group and the existence of a conditional expectation on $C^*(P)$ satisfying the conclusion in Lemma~\ref{fce}?

Furthermore, it would also be interesting to study the question for which subsemigroups of groups the left regular representation $\lambda: C^*_s(P) \to C^*_r(P)$ is an isomorphism. This is a weaker requirement than left amenability of $P$. Indeed, we have seen in Section~\ref{amenability} that the difference between the statements \an{$\lambda: C^*_s(P) \to C^*_r(P)$ is an isomorphism} and \an{$P$ is left amenable} is precisely given by the property of left reversibility. In this context, A. Nica has studied the example $P = \Nz^{*n}$, the $n$-fold free product of $\Nz$. He has shown in \cite{Ni}, Section~5 that although this semigroup is not left amenable, its left regular representation $\lambda: C^*(\Nz^{*n}) \to C^*_r(\Nz^{*n})$ is an isomorphism. So, the following question remains open: How can we characterize those semigroups which are not left amenable but still satisfy the condition that their left regular representations are isomorphisms?

Finally, let us come back to the construction of semigroup C*-algebras due to G. Murphy in \cite{Mur2} and \cite{Mur3} mentioned in the introduction. One could say that G. Murphy's construction leads to very complicated or even not tractable C*-algebras because the general theory of isometric semigroup representations is extremely complex. If we compare his construction with ours, then we see that G. Murphy's C*-algebras encode all isometric representations of the corresponding semigroups whereas representations of our C*-algebras correspond to rather special isometric representations because of the extra relations we have built into our construction. At the same time, these extra relations lead to a close relationship between our semigroup C*-algebras and the semigroups themselves in the context of amenability. Such a close relationship does not exist for G. Murphy's construction. For example, his semigroup C*-algebra of the semigroup $\Nz \times \Nz$ is by definition the universal C*-algebra generated by two commuting isometries. But this C*-algebra is not nuclear by \cite{Mur4}, Theorem~6.2. Such phenomena cannot occur in our setting by Corollary~\ref{C*(abelian)nuclear}.

\end{document}